\documentclass[12pt]{article}\def\today{March 20, 2003}
\usepackage{amssymb}
\usepackage{amsmath}
\usepackage{theorem}
\usepackage{latexsym}
\theoremheaderfont{\bf} 
\theorembodyfont{\sl}
\newtheorem{theo}{Theorem}[section]
{\theorembodyfont{\rm} \newtheorem{defi}[theo]{Definition}}
{\theorembodyfont{\rm} \newtheorem{exa}[theo]{Example}}
{\theorembodyfont{\rm} \newtheorem{rem}[theo]{Remark}}
\newtheorem{prop}[theo]{Proposition}
\newtheorem{cor}[theo]{Corollary}
\newtheorem{lemma}[theo]{Lemma}

\newtheorem{conjecture}[theo]{Conjecture}
\newenvironment{proof}{{\sc Proof:}}{\mbox{}\hfill$\Box$\par}
\newcommand{\eqr}[1]{~\mbox{$(${\rm \ref{#1}}$)$}}

\renewcommand{\theequation}{\thesection.\arabic{equation}}

\newcommand{\junk}[1]{}

\newcommand{\N}{{\mathbb N}}
\newcommand{\F}{{\mathbb F}}

\newcommand{\Z}{{\mathbb Z}}
\newcommand{\C}{{\mathcal C}}
\renewcommand{\H}{{\mathcal H}}
\newcommand{\wt}{{\rm wt}}
\newcommand{\rank}{{\rm rank}\,}
\newcommand{\spann}{\mbox{\rm span}}
\newcommand{\dfree}{\mbox{$d_{\mbox{\rm\tiny free}}$}}
\newcommand{\delay}[1]{\mbox{$\overleftarrow{#1}$}}
\newcommand{\T}{\mbox{$\!^{\sf T}$}}
\newcommand{\floor}[1]{\mbox{$\lfloor{#1}\rfloor$}}
\newcommand{\ceiling}[1]{\mbox{$\lceil{#1}\rceil$}}
\newcommand{\FDlaurent}{\mbox{$\F(\!(D)\!)$}}

\newcounter{abc}

\newenvironment{alphalist}{\begin{list}{(\alph{abc})\hfill}{\usecounter{abc}
     \topsep.5ex \labelwidth.6cm \leftmargin.7cm \labelsep.1cm
     \rightmargin0cm \parsep0ex \itemsep.6ex
     \partopsep1.6ex}}{\end{list}}
\newenvironment{arabiclist}{\begin{list}{(\arabic{abc})\hfill}{\usecounter{abc}
     \topsep.5ex \labelwidth.6cm \leftmargin.7cm \labelsep.1cm
     \rightmargin0cm \parsep0ex \itemsep.6ex
     \partopsep1.6ex}}{\end{list}}
\newenvironment{algo}{\begin{list}{{\rmfamily\bf Step~\arabic{abc}:}\hfill}{\usecounter{abc}
     \topsep.5ex \labelwidth1.7cm \leftmargin1.7cm \labelsep0cm
     \rightmargin0cm \parsep0ex \itemsep.6ex
     \partopsep1.6ex}}{\end{list}}
\newcommand{\vier}[4]{\left[ \begin{array}{cc}
                   #1 & #2 \\ #3 & #4 \end{array} \right]}

\title{Strongly MDS Convolutional Codes \footnote{The research of
    this paper was presented at the 2002 IEEE International
    Symposium on Information Theory in Lausanne, Switzerland,
    June 30--July 5, 2002 and at the International Symposium on
    the Mathematical Theory of Networks and Systems (MTNS),
    University of Notre Dame, August, 12--16 2002.  The authors
    were supported in part by NSF grants DMS-00-72383 and CCR-02-05310.}}
\date{\today}
\author{
  Heide Gluesing-Luerssen\\
  {\small Department of Mathematics}\vspace{-2mm}\\
  {\small University of Oldenburg}\vspace{-2mm}\\
  {\small P.O. Box 2503}\vspace{-2mm}\\
  {\small D-26111 Oldenburg, Germany}\vspace{-2mm}\\
  {\small {\em e-mail:} gluesing@mathematik.uni-oldenburg.de}
  \and
  Joachim Rosenthal\\
  {\small Department of Mathematics}\vspace{-2mm}\\
  {\small University of Notre Dame}\vspace{-2mm}\\
  {\small Notre Dame, Indiana 46556-5683, USA}\vspace{-2mm}\\
  {\small {\em e-mail:} Rosenthal.1@nd.edu} \and
  Roxana Smarandache\\
  {\small Department of Mathematical Sciences}\vspace{-2mm}\\
  {\small San Diego State University}\vspace{-2mm}\\
  {\small 5500 Campanile Dr.}\vspace{-2mm}\\
  {\small San Diego, CA 92182-7720, USA}\vspace{-2mm}\\
  {\small {\em e-mail:} rsmarand@sciences.sdsu.edu }}

\begin{document}
\maketitle
\begin{abstract}
  MDS convolutional codes have the property that their free
  distance is maximal among all codes of the same rate and the
  same degree. In this paper we introduce a class of MDS
  convolutional codes whose column distances reach the
  generalized Singleton bound at the earliest possible instant.
  We call these codes strongly MDS convolutional codes. It is
  shown that these codes can decode a maximum number of errors
  per time interval when compared with other convolutional codes
  of the same rate and degree. These codes have also a maximum or  near
  maximum distance profile. A code has a maximum distance profile
  if and only if the dual code has this property. \bigskip

\end{abstract}

\newpage

\section{Introduction}

In comparison to the literature on linear block codes there exist
only relatively few algebraic constructions of
convolutional codes having some good designed distance. There are
even fewer algebraic decoding algorithms which are capable of
exploiting the algebraic structure of the code.

Convolutional codes are typically decoded via the Viterbi
algorithm which has the advantage that soft information can be
processed. This algorithm has however the disadvantage that it is
too complex for codes with large degree or large memory or when
the block length is large. The algorithm is also not practical
for convolutional codes defined over large alphabets. There are
some alternative sub-optimal algorithms such as sequential
decoding and feedback decoding.  All these algorithms do not in general
exploit the algebraic structure of the convolutional code.

In applications where codes over large alphabets are required the
codes of choice are linear block codes with large distance such as
Reed-Solomon codes and more general algebraic geometric
codes. These codes can be algebraically decoded using e.g. the
Berlekamp-Massey algorithm or some of its generalizations.

In this paper we introduce a new class of convolutional codes
which we call {\em strongly MDS convolutional codes}. These codes
are particularly suited for applications where large alphabets
are involved. The free distance of these codes reaches the
generalized Singleton bound. This is the maximal possible
distance a convolutional code of a certain rate and degree can
have. The number of errors that strongly MDS convolutional codes can
correct per time interval is in a certain sense maximal as well.
We will make this precise in Section~\ref{Sec6}.

Let $\F$ be any finite field and denote by $\F[D]$ and
$\FDlaurent$ the polynomial ring respectively the field of all
formal Laurent series over $\F$, i.~e.
\[
  \F[D]=\Big\{\sum_{j=0}^L a_jD^j\,\Big|\, L\in\N_0, a_j\in\F\Big\}\text{ and }
  \FDlaurent=\Big\{\sum_{j=l}^{\infty} a_j D^j\,\Big|\, l\in\Z, a_j\in\F\Big\}.
\]
For $v=\sum_{j=l}^{\infty}v_jD^j\in\FDlaurent^n\backslash\{0\}$ we define
$\delay{v}:=\min\{j\in\Z\mid v_j\not=0\}$ to be the {\em delay\/}
of the sequence~$v$, that is the time instant, at which the sequence
actually starts.
We put $\delay{0}:=\infty$.

Let $G\in\F[D]^{k\times n}$ be a $k\times n$ polynomial matrix of
rank $k$.  We define a convolutional code of rate $k/n$ as the
set
\begin{equation}\label{e-convcode}
   \C:=\{uG\mid u\in\FDlaurent^k\}\subseteq\FDlaurent^n
\end{equation}
and say that $G$ is a generator matrix of the code $\C$.
Two generators of $\C$ differ only by a nonsingular left transformation over
$\FDlaurent$.
It is well-known that we can assume $G$ to be basic and minimal
in the following sense.

\begin{defi}[see~\cite{fo70}]\label{D-basic.minimal}
  A polynomial generator matrix $G\in\F[D]^{k\times n}$ is called
  {\em basic} if it has a polynomial right inverse (equivalently,
  if the $k\times k$-minors are coprime in $\F[D]$).  It is
  called {\em minimal} if $\sum_{i=1}^k\nu_i$, where $\nu_i$
  denotes the $i$th row degree of $G$, attains the minimal value
  among all generator matrices of $\C$.
\end{defi}

Two basic generator matrices differ only by a unimodular left
transformation over $\F[D]$.  If~$G$ is a minimal basic encoder
one defines the {\em degree}~\cite{mc98} of $\C$ as the number
$\delta:=\sum_{i=1}^k\nu_i$.  In the literature the
degree~$\delta$ is sometimes also called the {\em total
  memory}~\cite{li83} or the {\em overall constraint
  length}~\cite{jo99} or the {\em complexity}~\cite{pi88} of the
minimal basic generator matrix $G(D)$.  We like to use the term
degree as it corresponds to the term {\em
McMillan degree\/} used in systems theory~\cite{fo75,ro96a1,ro99a}.
We also wish to point out that in algebraic geometry the degree corresponds
to the degree of an associated vector bundle
(i.e. quotient sheaf), see~\cite{lo01,ro01,ro99a1} for more
details.

Since the degree depends only on the code itself, but not on the
specific choice of the generator matrix $G$, we will
call~$\delta$ the degree of the code~$\C$.  Recall also from
Forney~\cite{fo75} that the set $\{\nu_1,\ldots,\nu_k\}$ of row
degrees is the same for all minimal basic encoders of $\C$.
Because of this reason McEliece~\cite{mc98} calls these indices
the {\em Forney indices\/} of the code $\C$.  As a consequence, also the
number $\nu:=\max\{\nu_1,\ldots,\nu_k\}$ depends only on the code
$\C$ itself and is usually called the {\em memory\/} of the code.
In the sequel we will adopt the notation used by
McEliece~\cite[p.~1082]{mc98} and call a convolutional code of
rate $k/n$ and degree $\delta$ an $(n,k,\delta)$-code.  Every
$(n,k,\delta)$-code~$\C$ can also be represented in terms of a
parity check matrix, i.~e. a matrix $H\in\FDlaurent^{(n-k)\times
  n}$ such that
\[
    \C=\{v\in\FDlaurent^n\mid vH\T=0\}.
\]
It is clear that we can choose $H$ to be polynomial, thus
$H\in\F[D]^{(n-k)\times n}$, and basic.
Notice also that $GH\T=0$ for any generator matrix~$G$ of~$\C$.

For a vector $v\in \F^n$, we define its weight
$\wt(v)$ as the number of all its nonzero components.
For $v=\sum_{j=l}^{\infty} v_jD^j\in\FDlaurent^n$ we define
\[
  \wt(v):=\sum_{j=l}^{\infty}\wt(v_j)\in\N_0\cup\{\infty\}.
\]
Finally, the {\em free distance} of the convolutional code
$\C\subset \FDlaurent^n$ is defined through
\begin{equation}\label{e-dfree}
   \dfree:=\min\{\wt(v)\mid v\in \C, v\neq 0\}.
\end{equation}
It is an easy, but crucial observation that a basic generator matrix $G$
yields a non-catastrophic and delay-free encoder, i.~e., if
$v=uG\in\FDlaurent^n$ for some $u\in\FDlaurent^k$, then
\[
   \wt(v)\text{ finite }\Longrightarrow \wt(u)\text{ finite}
\]
and
\[
  \delay{v}=\delay{u}.
\]
Therefore, in case we are given a basic generator matrix $G$, the free
distance can also be obtained as
\begin{align*}
  \dfree&=\min\big\{\wt(v)\,\big|\, v=uG\text{ for some }
               u\in\F[D]^k\backslash\{0\}\big\}\\
        &=\min\big\{\wt(v)\,\big|\, v=uG\text{ for some }
               u\in\F[D]^k\backslash\{0\},u_0\not=0\big\}.
\end{align*}

An $(n,k,\delta)$ convolutional code is called MDS if its free distance
is maximal among all rate $k/n$ convolutional codes of degree
$\delta$, i.e. an  $(n,k,\delta)$ convolutional code is MDS if
the free distance achieves the generalized Singleton bound~\cite{ro99a1}:
$$
\dfree= (n-k)\Big(\Big\lfloor\frac{\delta}{k} \Big\rfloor +1
\Big)+\delta +1.
$$
The concept of MDS convolutional codes was introduced by the
authors in~\cite{ro99a1,sm01a}.  Strongly MDS codes are going to
be a subclass of MDS codes which have a remarkable decoding
capability.

The paper is structured as follows: In Section~\ref{Sec2} we
review notions from convolutional coding theory such as the
column distances, the generalized Singleton bound and we
introduce the important concepts for this paper, namely the
property of being strongly MDS and having a maximum distance
profile.  In Section~\ref{Sec3} we show the existence of strongly
MDS codes in the situation when the rate is $(n-1)/n$. In order
to do so we introduce the interesting concept of a {\em
  superregular matrix} which might be of independent interest. In
Section~\ref{Sec4} we illustrate the concepts through a series of
examples.  In Section~\ref{Sec5} we investigate to what extend
properties of MDS, strongly MDS and maximum distance profile
carry over to the dual code. The main result of this section
states that a code has a maximum distance profile if and only if
its dual has this property. This allows us then to show that for
certain specific parameters a code is strongly MDS if and only if
its dual is strongly MDS. Finally in Section~\ref{Sec6} we show
how strongly MDS convolutional codes can be decoded via feedback
decoding. It turns out that the number of errors which can be
decoded per time interval compares well to a maximum distance
separable block code.

\section{Strongly MDS Codes and Codes with Maximum Distance Profile}
\setcounter{equation}{0}     \label{Sec2}
In this section we will recall the column distances of a convolutional
code and their relation to the free distance.
After showing some upper bounds for these distances we will introduce
the notion of strongly MDS codes.
It describes codes, for which the column distances attain their maximum
value.

Throughout this section let $\C\subseteq\FDlaurent^n$ be an
$(n,k,\delta)$-code with basic generator matrix
\begin{equation}\label{e-G}
   G=\sum_{j=0}^\nu G_jD^j\in\F[D]^{k\times n},\ G_j\in\F^{k\times n},
   G_{\nu}\not=0
\end{equation}
and basic parity check matrix
\begin{equation}\label{e-H}
  H=\sum_{j=0}^{\mu}H_jD^j\in\F[D]^{(n-k)\times n},\
  H_j\in\F^{(n-k)\times n}, H_{\mu}\not=0.
\end{equation}
Notice that $\nu$ is the memory of the code.
For every $j\in\N_0$ we define the truncated sliding generator and parity check matrices
\begin{equation}\label{e-Gcj}
\begin{array}{rcl}
       G^c_j &= &\begin{bmatrix}
                 G_0& G_1& \ldots& G_j\\
                    & G_0& \ldots& G_{j-1}\\
                    &    & \ddots& \vdots\\
                    &    &       & G_0
              \end{bmatrix}\in\F^{(j+1)k\times(j+1)n},\\
       H^c_j: &=& \begin{bmatrix}
                 H_0&    &   &      \\
                 H_1& H_0&   &      \\
                \vdots&\vdots&\ddots&\\
                 H_j&H_{j-1} &\ldots& H_0
              \end{bmatrix}\in\F^{(j+1)(n-k)\times(j+1)n},
\end{array}
\end{equation}
where we let $G_j=0$ (resp.\ $H_j=0$) whenever $j>\nu$ (resp.\ $j>\mu$),
see also \cite[p.~110]{jo99}.
The identity $GH\T=0$ immediately implies $G^c_j(H^c_j)\T=0$ for all
$j\in\N_0$.
Since~$G$ and~$H$ are both basic, the matrices~$G^c_j$ and~$H^c_j$ both
have full rank and therefore we even have
\begin{equation}\label{e-im.ker}
   \{uG^c_j\mid u\in\F^{(j+1)k}\}=\{v\in\F^{(j+1)n}\mid v (H^c_j)\T=0\}
   \text{ for all }j\in\N_0.
\end{equation}
The relevance of these matrices rests on the fact that they single out
codeword sequences of length~$j$ in the following sense.

\begin{rem}\label{R-GcjHcj}
For $v:=\sum_{j=l}^{\infty} v_jD^j\in\FDlaurent^n$ and $m,M\in\Z$ with
$m\leq M$ define
\[
   v_{[m,M]}:=(v_m,v_{m+1},\ldots,v_M)\in\F^{(M-m+1)n}.
\]
Then we have the following:
\begin{alphalist}
\item If $v=uG$ for some $u\in\FDlaurent^k$ with $\delay{u}\geq l$, then
      $v_{[l,l+j]}=u_{[l,l+j]}G^c_j$ and $v_{[l,l+j]}(H^c_j)\T=0$ for all
      $j\in\N_0$.
\item If $\hat{v}=\hat{u}G^c_j\in\F^{(j+1)n}$ for some
      $j\in\N_0$ and $\hat{u}\in\F^{(j+1)k}$, then there exists $v\in\C$
      such that $\delay{v}\geq0$ and $v_{[0,j]}=\hat{v}$.
\item For all $j\in\N_0$ we have
      $\{v_{[0,j]}\mid v\in\C,\delay{v}=0\}
       =\{\hat{v}=(\hat{v}_0,\ldots,\hat{v}_j)\in\F^{(j+1)n}\mid
        \hat{v}(H^c_j)\T=0,\,\hat{v}_0\not=0\}$.
\end{alphalist}
Part~(a) follows easily by equating like powers of~$D$ in the equation
$v=uG$ and by use of\eqr{e-im.ker};
(b)~is obvious by taking $v=uG$ with $u=\sum_{i=0}^j \hat{u}_iD^i$, where
$\hat{u}=(\hat{u}_0,\ldots,\hat{u}_j)$; (c) is a consequence
of~(a) and\eqr{e-im.ker}.
\end{rem}

Following \cite[pp.~110]{jo99}, the $j$th {\em column distance of the code\/}~$\C$
is defined to be
\begin{equation}\label{e-dcj2}
  d^c_j:=\min\big\{\wt(v_{[0,j]})\,\big|\, v\in\C,
                     \delay{v}=0\}.
\end{equation}
Using the remark above we obtain the alternative identities
\begin{align}
    d^c_j
       &=\min\big\{\wt(v_{[0,j]})\,\big|\,
                v=uG, u\in\F[D]^k, u_0\not=0\big\}\nonumber\\[.6ex]
       &=\min\big\{\wt(v_{[l,l+j]})\,\big|\,
                   v=uG,\,u\in\FDlaurent^k,\,\delay{u}=l\big\}\nonumber\\[.6ex]
       &=\min\big\{\wt\big((u_0,\ldots,u_j)G^c_j\big)\,\big|\,
                  u_i\in\F^k,\,u_0\not=0\big\}\label{e-dcj}\\[.6ex]
       &=\min\{\wt(\hat{v})\mid
       \hat{v}=(\hat{v}_0,\ldots,\hat{v}_j)\in\F^{(j+1)n},\,
            \hat{v}(H^c_j)\T=0,\,\hat{v}_0\not=0\}.\label{e-dcj3}
\end{align}
Obviously, $d^c_j\leq\dfree$ for all $j\in\N_0$ and one even
has~\cite[pp.~113]{jo99}
\begin{equation}\label{e-dist.inequ}
  d^c_0\leq d^c_1\leq d^c_2\ldots \text{ and }\lim_{j\rightarrow\infty}d^c_j=\dfree.
\end{equation}
The $(\nu +1)$-tuple  of numbers $(d^c_0, d^c_1,
d^c_1,\ldots,d^c_\nu)$, where $\nu$ is the memory, is
is called the {\em column distance profile} of the
code~\cite[p. 112]{jo99}.

Equation\eqr{e-dcj3} immediately implies

\begin{prop}\label{P-dcj}
Let $d\in\N$. Then the following properties are equivalent.
\begin{alphalist}
\item $d^c_j=d$;
\item none of the first~$n$ columns of~$H^c_j$ is contained in the span of any other
       $d-2$ columns and one of the first~$n$ columns of~$H^c_j$ is in the span of
       some other $d-1$ columns of that matrix.
\end{alphalist}
\end{prop}
We leave it to the reader to verify the equivalence of the statements.

\begin{prop}\label{P-dcj.bound}
For every $j\in\N_0$ we have
\[
    d^c_j\leq(n-k)(j+1)+1.
\]
\end{prop}

\begin{proof}
  Consider the sliding parity check matrix $H^c_j$ introduced
  in\eqr{e-Gcj}. The set of vectors
  $\hat{v}=(\hat{v}_0,\ldots,\hat{v}_j)\in\F^{(j+1)n},\,
  \hat{v}(H^c_j)\T=0,\,\hat{v}_0\not=0$ forms a nonlinear subset
  of the linear block code defined by the left kernel of $H^c_j$.
  Since $H^c_j\in\F^{(n-k)(j+1)\times n(j+1)}$ any vector in the
  left kernel has weight at most $(n-k)(j+1)+1$ by the usual
  Singleton bound for block codes and this establishes the claim.
\end{proof}

The column distances give information about the error-correcting
capabilities of the code.
Precisely, $d^c_j$ determines the error-correcting capability of a decoder that
estimates the message symbol $u_0$ based on the received symbols
$v_{[0,j]}$, see also \cite[p.~111]{jo99}.
Therefore, a good performance for sequential decoding requires the column
distances as big as possible.
The next proposition shows that maximality of $d^c_j$ implies maximality
of the preceding column distances.

\begin{cor}\label{C-optdistances}
If $d^c_j=(n-k)(j+1)+1$ for some $j\in\N_0$, then
$d^c_i=(n-k)(i+1)+1$ for all $i\leq j$.
\end{cor}

\begin{proof}
It suffices to prove the assertion for $i=j-1$.
In order to do so notice that
\[
    H^c_j=\left[\!\begin{array}{c|c}
                & 0 \\
            H^c_{j-1}& \vdots \\
                     & 0\\ \hline
            H_j\:H_{j-1}\:\cdots\:H_1&H_0\rule{0cm}{2.5ex}
          \end{array}\!\right]
\]
and assume that one of the first~$n$ columns of~$H^c_{j-1}$ is in the span of some other
$(n-k)j-1$ columns.
Then $\rank H_0=n-k$ implies that one of the first~$n$ columns of~$H^c_j$ is in the span of
some other $(n-k)j-1+n-k=(n-k)(j+1)-1$ columns of $H^c_j$.
But this is a contradiction to the optimality of $d^c_j$ by
Proposition~\ref{P-dcj}.
\end{proof}

The Singleton-bound for block codes has been generalized to
convolutional codes in~\cite{ro99a1}. Therein the following has been
shown.

\begin{theo}\label{T-MDS}
The free distance of an $(n,k,\delta)$-code satisfies
\begin{equation}                                \label{G-Singleton}
   \dfree\leq (n-k)\Big(\Big\lfloor\frac{\delta}{k}\Big\rfloor+1\Big)+\delta+1.
\end{equation}
\end{theo}
The number appearing on the right in\eqr{G-Singleton} is called the {\em
generalized Singleton bound}.
The code is called an MDS code if it satisfies
$\dfree=(n-k)\big(\floor{\frac{\delta}{k}}+1\big)+\delta+1$.
It has been shown in~\cite{ro99a1} that
for every set of parameters $(n,k,\delta)$ and every
prime number~$p$ there exists a suitably large finite field~$\F$
of characteristic~$p$ and an MDS code with parameters $(n,k,\delta)$ over $\F$.

The generalized Singleton bound reduces to the usual Singleton
bound $n-k+1$ when $\delta=0$, the block code situation.

The proof of the existence of MDS codes given in~\cite{ro99a1} is
based on techniques from algebraic geometry and is
non-constructive.  In~\cite{sm01a} a construction of MDS codes
with parameters $(n,k,\delta)$ was given for suitably large
fields of characteristic coprime with~$n$.

In the sequel we will strengthen the MDS property by requiring that the generalized
Singleton bound is
attained by the earliest column distance possible.
This will lead to the notion of a strongly MDS code.

\begin{prop}\label{P-bound.for.j}
Suppose~$\C$ be an MDS code with parameters $(n,k,\delta)$, column
distances $d^c_j,\,j\in\N_0,$ and free distance $\dfree$.
Let $M:=\min\{j\in\N_0\mid d^c_j=\dfree\}$. Then
\[
   M\geq\Big\lfloor\frac{\delta}{k}\Big\rfloor+\Big\lceil\frac{\delta}{n-k}\Big\rceil.
\]
\end{prop}

\begin{proof}
From Proposition~\ref{P-dcj.bound} we get
\begin{equation}\label{e-dcM}
    \dfree=(n-k)\Big(\Big\lfloor\frac{\delta}{k}\Big\rfloor+1\Big)+\delta+1=d^c_M\leq
  (n-k)(M+1)+1.
\end{equation}
This yields the assertion.
\end{proof}

The proof also shows that in the case
$j>\floor{\frac{\delta}{k}}+\ceiling{\frac{\delta}{n-k}}$ the column
distance~$d^c_j$ never attains the upper bound
$(n-k)(j+1)+1$ of Proposition~\ref{P-dcj.bound}, see also\eqr{e-dist.inequ}.

\begin{defi}\label{D-strongly.MDS}
  An $(n,k,\delta)$-code with column distances
  $d^c_j,\,j\in\N_0$, is called {\em strongly} MDS, if
\[
   d^c_M=(n-k)\Big(\Big\lfloor\frac{\delta}{k}\Big\rfloor+1\Big)+\delta+1
   \text{ for } M=\Big\lfloor\frac{\delta}{k}\Big\rfloor+\Big\lceil\frac{\delta}{n-k}\Big\rceil.
\]
\end{defi}
Because of\eqr{e-dist.inequ} the strong MDS property implies the
MDS property.

\begin{rem}\label{R-strMDS}
In the case where $(n-k)\mid\delta$, the strong
MDS property implies that $d^c_M$ attains the upper bound
$(n-k)(M+1)+1$, see Proposition~\ref{P-dcj}.
Hence Corollary~\ref{C-optdistances} implies that in this case {\em all\/}
column distances attain their optimal value.
\\
If $(n-k)\nmid\delta$, we always have $d^c_M<(n-k)(M+1)+1$ as can be seen
from\eqr{e-dcM}.
\end{rem}

Even when $(n-k)\nmid\delta$ it is very desirable that the column
distance profile $ d^c_0, d^c_1, d^c_1,\ldots$ has the maximum
possible increase at each step. This motivates the following
definition.
\begin{defi}           \label{max-Profile}
Let
\begin{equation}                                         \label{e-L}
L:=\Big\lfloor\frac{\delta}{k}\Big\rfloor+\Big\lfloor\frac{\delta}{n-k}\Big\rfloor.
\end{equation}
An $(n,k,\delta)$-code with column distances $d^c_j,\,j\in\N_0$, is said
to have a {\em maximum distance profile} if
\[
   d^c_j=(n-k)(j+1)+1,\ \mbox{ for }\ j=1,\ldots,L.
\]
\end{defi}
Using the notation of Definition~\ref{D-strongly.MDS} we have
\begin{equation}\label{e-LM}
    L=\left\{\begin{array}{ll} M &\text{if } (n-k)\mid\delta\\[1ex]
                       M-1&\text{otherwise.}
      \end{array}\right.
\end{equation}
An immediate consequence of Corollary~\ref{C-optdistances} is
\begin{lemma}\label{L-MDP}
   An $(n,k,\delta)$-code has a maximum distance profile if and
   only if the $L$th column distance satisfies
  \[
   d^c_L=(n-k)(L+1)+1.
\]
\end{lemma}
As a consequence we obtain that if $n-k$ divides~$\delta$ then an $(n,k,\delta)$-code has
maximum distance profile if and only if it is strongly MDS since
$(n-k)\big(\big\lfloor\frac{\delta}{k}\big\rfloor+\frac{\delta}{n-k}+1\big)+1=
 (n-k)(\big\lfloor\frac{\delta}{k}\big\rfloor+1\big)+\delta+1$.

\begin{rem}
  The concept is clearly related to the notion of {\em optimum
    distance profile} (ODP), see~\cite[p.~112]{jo99}. For ODP it
  is required that the column distances are maximal up to the
  memory $\nu$. Hence if $\nu\leq L$ then a code with maximum
  distance profile is always ODP. In general one expects a good
  code to have generic Forney indices, i.e. the indices attain
  only the two values $\lceil \frac{\delta}{k}\rceil$ and
  $\lfloor \frac{\delta}{k}\rfloor$.
  McEliece~\cite[Corollary~4.3]{mc98} calls such codes {\em
    compact codes}. It has been shown in~\cite{ro99a1} that an
  MDS code has always generic indices.  Of course if the indices
  are generic then $\nu=\lceil \frac{\delta}{k}\rceil$ and thus
  $\nu\leq L+1$.

  The notion of ODP seems also to be dependent on the base field
  which is usually assumed to be the binary field. A code with
  maximum distance profile does in general not exist over the
  binary field and it can only exist for sufficiently large base
  fields. This is similar to the situation of MDS block codes.
  Such codes are known to exist as soon as the field size of $\F$
  is larger than the block length $n$.
\end{rem}

One of the main results of Section~\ref{Sec5} will show that a
convolutional code has a maximum distance profile if and only if
its dual has this property. The following algebraic criterion
which characterizes codes having a maximum distance profile will
be very useful.

\begin{theo}                        \label{Alg-criterion}
 Let $G=\sum_{j=0}^\nu G_jD^j$ be the generator matrix of an
 $(n,k,\delta)$-code. Let $L$ be defined as in\eqr{e-L} and let
 \begin{equation}
\begin{array}{rcl}
       G^c_L &= &\begin{bmatrix}
                 G_0& G_1& \ldots& G_L\\
                    & G_0& \ldots& G_{L-1}\\
                    &    & \ddots& \vdots\\
                    &    &       & G_0
              \end{bmatrix}\in\F^{(L+1)k\times(L+1)n}.
\end{array}
\end{equation}
Then $G$ represents a maximum distance profile code if and only
if every $(L+1)k\times(L+1)k$ full-size minor formed from the
columns with indices $1\leq j_1< \ldots < j_{(L+1)k}$, where
$j_{sk+1}>sn$ for $s=1,\ldots,L$, is nonzero.
\end{theo}

\begin{proof}
  Assume there are indices $1\leq j_1< \ldots < j_{(L+1)k}$
  satisfying $j_{sk+1}>sn$ for $s=1,\ldots,L$ whose corresponding
  minor is zero. It follows that there is a vector
  $u=(u_0,\ldots,u_L)$ such that $uG^c_L$ has zero coordinates at
  positions $j_1,\ldots , j_{(L+1)k}$. Let $\ell :=\min \{ i\mid
    u_i\neq 0\}$. Consider the vector
  $$
   \left(  u_\ell ,\ldots,u_L \right) G^c_{L-\ell}\in\F^{(L-\ell +1)n}.
  $$
  The weight of this vector is at most $(L-\ell +1)(n-k)$ as
  there are at least  $(L-\ell +1)k$ coordinates zero. It follows
  from\eqr{e-dcj} that $d^c_{L-\ell}\leq (L-\ell+1)(n-k)$ and by
  Corollary~\ref{C-optdistances} the code has not a maximum
  distance profile.\medskip

  Vice versa assume that $\C$ has not a maximum distance profile.
  Let $m :=\min \{ i\mid d_i^c\leq (n-k)(i+1)\}$. It follows that
  there is a vector $u=(u_0,\ldots,u_m)$, $u_0\neq 0$ such that
  $uG^c_m$ has at least $k(m+1)$ zeros. As a submatrix inside
  $G^c_L$ we select the columns corresponding to the first
  $k(m+1)$ positions where $uG^c_m$ has a zero and we augment it
  by the last $k(L-m)$ columns of $G^c_L$.
  We call the indices of the selected columns $j_1,\ldots,j_{(L+1)k}$.
  This gives an
  $(L+1)k\times(L+1)k$ full-size minor and we claim that this
  minor is zero and that the indices $j_1,\ldots , j_{(L+1)k}$
  satisfy $j_{sk+1}>sn$ for $s=1,\ldots,L$. In order to prove the
  latter note that $d_i^c= (n-k)(i+1)+1$ for $i=0,\ldots,m-1$. It
  therefore follows that $(u_0,\ldots,u_i)G^c_i$ has at most
  $k(i+1)-1$ zeros for $i=0,\ldots,m-1$. In particular  $j_{sk+1}>sn$
  for $s=1,\ldots,m$. Clearly it is also true for $s=m+1,\ldots,
  L$. It remains to be shown that the minor is zero. For this
  note that the selected matrix has the form $\vier{A}{B}{0}{C}$
  where $A$ is an $(m+1)k\times(m+1)k$ submatrix of $G^c_m$ which
  is singular by construction. The fullsize minor is therefore zero.
\end{proof}

\section{Existence of Strongly MDS \protect{$(n,n-1,\delta)$}-Codes}
\setcounter{equation}{0}                           \label{Sec3}

During his investigation of algebraic decoding of convolutional
codes B. Allen conjectured in his dissertation~\cite{al99t} the
existence of strongly MDS convolutional codes in the situation
when $k=1$ and $n=2$.  In this section we will show the existence
of strongly MDS codes with parameters $(n,n-1,\delta)$.
It follows from Equation\eqr{e-LM} and Lemma~\ref{L-MDP} that these codes
also have maximum distance profile.
By Theorem~\ref{T-MDS} the generalized Singleton bound for these parameters
is given by $\floor{\frac{\delta}{n-1}}+\delta+2$.  Thus,
Definition~\ref{D-strongly.MDS} yields that we have to find an
$(n,n-1,\delta)$-code such that $d_M^c=M+2$, where
$M=\floor{\frac{\delta}{n-1}}+\delta$.  In order to do so, let
\begin{equation}\label{e-Hn}
    H=[a_1,\ldots,a_n]\in\F[D]^{n}, \text{ where }
    a_i=\sum_{j=0}^{\delta}a_{ij}D^j\in\F[D],
\end{equation}
be a basic parity check matrix of the desired code.
Without loss of generality we may assume $a_{10}=1$.
The strong MDS property can now be expressed as follows.

\begin{theo}\label{T-strongMDS}
Let $H\in\F[D]^{n}$ be as in\eqr{e-Hn}, let $a_{10}=1$ and define
$\C:=\{v\in\FDlaurent^n\mid vH\T=0\}$ be the code with parity check matrix $H$.
Furthermore, for $i=2,\ldots,n$ let
\begin{equation}\label{e-Laurent2}
   \frac{a_i}{a_1}=\sum_{j=0}^{\infty}h_{ji}D^j\in\FDlaurent
\end{equation}
be the Laurent expansion of $\frac{a_i}{a_1}\in\F(D)$ and for
$M=\lfloor\frac{\delta}{n-1}\rfloor+\delta$ define
\begin{eqnarray}
  \hat{H}:=& \label{e-hatH1}
    \left[\!\begin{array}{cccc|cccc}
      1&      &      & &  h_{02}\ \cdots\  h_{0n}&          &   &  \\
       &\ddots&      & &  h_{12}\ \cdots\ h_{1n}&h_{02}\quad \cdots\quad h_{0n}& & \\
       &      &\ddots& & \vdots\qquad\quad\vdots&\vdots\qquad\qquad\quad\vdots&\ddots &  \\
       &      &      &1& h_{M2}\cdots\ h_{Mn}&h_{M-1,2}\hfill\cdots\ h_{M-1,n}
                                         &\cdots\cdots&h_{02}\ \cdots\ h_{0n}
  \end{array}\!\right]\\[2.5ex]
     =:&\label{e-hatH2}
       [e_1,\ldots,e_{M+1},H_{12},\ldots,H_{1n},\ldots,H_{M+1,2},\ldots,H_{M+1,n}]
       \in\F^{(M+1)\times(M+1)n},
\end{eqnarray}
where $e_i$ denotes the $i$-th standard basis vector.
We call $\hat{H}$ the $M$-th systematic sliding parity check matrix of~$\C$.
The following conditions are equivalent:
\begin{alphalist}
\item $\C$ is strongly MDS, i.~e. $d^c_M=M+2$,
\item none of the columns $H_{12},\ldots, H_{1n}$ of~$\hat{H}$
      is contained in the span of any other~$M$ columns of~$\hat{H}$.
\end{alphalist}
Notice that~(b) automatically implies that $h_{ji}\not=0$ for all~$i$ and~$j$ and that
also the first column $e_1$ is not in the span of any other~$M$ columns.
\end{theo}

\begin{proof}
After a column permutation the sliding parity check matrix
$H^c_M\in\F^{(M+1)\times(M+1)n}$ of $\C$ has the form
\[
   H':=\left[\!\begin{array}{cccc|cccc}
       1&      &   &   &   a_{20}\ \cdots\  a_{n0}&          &   &  \\
       a_{11}&\ddots&   &   &   a_{21}\ \cdots\ a_{n1}&a_{20}\quad \cdots\quad a_{n0}& & \\
       \vdots&\ddots&\ddots& &   \vdots\qquad\quad\vdots&\vdots\qquad\qquad\quad\vdots&\ddots &  \\
  a_{1M}&\ldots&a_{11}&1& a_{2M}\cdots\ a_{nM}&a_{2,M-1}\hfill\cdots\ a_{n,M-1}
                                         &\cdots\cdots&a_{20}\ \cdots\ a_{n0}
  \end{array}\!\right].
\]
It is straightforward to see that left multiplication of $H'$ by the inverse of the
first block leads to the matrix $\hat{H}$ of\eqr{e-hatH1}.
After these transformations Proposition~\ref{P-dcj} applied to the case $j=M$
and $d=M+2$ translates into the equivalence:
$\C$ is strongly MDS iff neither the first column~$e_1$ nor any of the columns
$H_{12},\ldots,H_{1n}$ is in the span of any other~$M$ columns of~$\hat{H}$.
But this in turn is equivalent to~(b) above.
\end{proof}

In order to establish the existence of strongly MDS codes we will proceed as
follows.
Firstly, we will establish the existence of a systematic sliding parity check
matrix~$\hat{H}$ as in\eqr{e-hatH1} with property~(b) of the theorem above.
Thereafter, we will show that there exist coprime polynomials $a_1,\ldots,a_n$ having
maximum degree equal to~$\delta$ such that
\[
   \frac{a_i}{a_1}=\sum_{j=0}^M h_{ji}D^j +\text{ higher powers},\ i=2,\ldots,n.
\]
Theorem~\ref{T-strongMDS} then yields that the code with parity check matrix
$H=[a_1,\ldots,a_n]$ is a strongly MDS $(n,n-1,\delta)$-code.

As for the first step, let us have a look at the special case of
$(2,1,\delta)$-codes. In this case $M=2\delta$ and the systematic sliding parity check matrix
in\eqr{e-hatH1} has the form
\begin{equation}\label{e-hatH}
    \hat{H}:=\left[\!\begin{array}{cccc|cccc}
                   1&  &      &  &h_0    & 0      & \cdots & 0 \\
                    & 1&      &  &h_1    & h_0   & \ddots & \vdots \\
                    &  &\ddots&  & \vdots & \ddots &\ddots & 0\\
                    &  &      &1 &h_{2\delta}& \cdots &   h_1&h_0\end{array}\!\right]
    =:[I_{2\delta+1}, T]
    \in\F^{(2\delta+1)\times(4\delta+2)}, \text{ where }
    h_j\in\F.
\end{equation}
As we will see, the existence of matrices~$T$ of any given size and the structure above
such that~$\hat{H}$ has the
column property of Theorem~\ref{T-strongMDS}(b) will be the main tool for the
existence of strongly MDS codes even of length $n>2$.
Therefore we will concentrate on these matrices first. The main point is to express
the column condition on~$\hat{H}$ in terms of the minors of~$T$.

\begin{defi}\label{D-minors}
Let $R$ be a ring.
For a matrix $T\in R^{n\times k}$ denote by
$T^{i_1,\ldots,i_r}_{j_1,\ldots,j_s}\in R^{r\times s}$ the $r\times s$-submatrix
obtained from~$T$ by picking the rows with indices $i_1,\ldots,i_r$
and the columns with indices $j_1,\ldots,j_s$.
\end{defi}

In the sequel the following property will play a crucial role.

\begin{defi}\label{D-superreg}
  Let $\F$ be field.  A lower triangular matrix $T\in\F^{n\times
    k}$ is said to be {\em superregular}\footnote{We adopt this
    notion from~\cite{ro89a2}, where it has been coined in a
    slightly different context.}, if
  $T^{i_1,\ldots,i_r}_{j_1,\ldots,j_r}$ is nonsingular for all
  $1\leq r\leq\min\{k,n\}$ and all indices $1\leq
  i_1<\ldots<i_r\leq n,\,1\leq j_1<\ldots<j_r\leq k$ which
  satisfy $j_{\nu}\leq i_{\nu}$ for $\nu=1,\ldots,r$.  We call
  the submatrices obtained by picking such indices the proper
  submatrices and their determinants the proper minors of~$T$.
\end{defi}

\begin{rem}\label{R-proper}
Observe that the proper submatrices are the only submatrices which can
possibly be nonsingular.
This can be seen as follows.
If $j_{\nu}>i_{\nu}$ for some $\nu$, then in
the submatrix $\hat{T}:=T^{i_1,\ldots,i_r}_{j_1,\ldots,j_r}$
the upper right block consisting of the first $\nu$ rows and the last $r-\nu+1$
columns is identically zero.
Hence the first $\nu$ rows of $\hat{T}$ can have at most rank $\nu-1$.
In other words, the improper submatrices of~$T$ are trivially singular.
For example, for $T=(h_{ij})$ we have
\[
    T^{1,2,5}_{1,3,4}=\begin{bmatrix} h_{11}&0&0\\ h_{21}&0&0\\
    h_{51}&h_{53}&h_{54}\end{bmatrix}.
\]
\end{rem}

Now we can establish the following.

\begin{theo}\label{T-superreg}
Let $\F$ be a field and~$T$ be a lower triangular Toeplitz matrix, i.~e.
\begin{equation}\label{e-HToeplitz}
    T=[T_1,\ldots,T_l]=
    \begin{bmatrix} h_0    & 0      & \cdots & 0 \\
                    h_1    & h_0    & \ddots & \vdots  \\
                    \vdots & \ddots &\ddots & 0\\
                    h_{l-1}    & \cdots &   h_1  & h_0\end{bmatrix}
    \in\F^{l\times l}.
\end{equation}
Furthermore, put
$\hat{H}:=[I_l,T]=[e_1,\ldots,e_l,T_1,\ldots,T_l]\in\F^{l\times2l}$.
Then the following are equivalent:
\begin{alphalist}
\item $T$ is superregular, i.e. all proper submatrices in the
  sense of Definition~\ref{D-superreg} are nonsingular.
\item Assume there are indices $1\leq i_1<\ldots<i_r\leq
  n,\,1< j_2<\ldots<j_r\leq k$. Then all proper submatrices
  of~$T$ of the form $T^{i_1,i_2,\ldots,i_r}_{1,j_2,\ldots,j_r}$
  are nonsingular,
\item $\wt\big(T_1+\sum_{j=1}^s\beta_j T_{m_j}\big)\geq l-s$ for all $1\leq s\leq l-1$, all
      $1<m_1<\ldots< m_s\leq l$ and all $\beta_1,\ldots,\beta_s\in\F$,
\item $T_1\not\in\spann\{T_{m_1},\ldots,T_{m_s},e_{l_1},\ldots,e_{l_t}\}$
      where $1< m_1<\ldots<m_s\leq l$ and $1\leq l_1<\ldots<l_t\leq l$ and
      $s+t\leq l-1$.
\item If $v\in\F^{2l}$ satisfies $v\hat{H}\T=0$ and $v_{l+1}\not=0$, then $\wt(v)\geq l+1$.
\item $e_1\not\in\spann\{T_{m_1},\ldots,T_{m_s},e_{l_1},\ldots,e_{l_t}\}$
      where $1\leq m_1<\ldots<m_s\leq l$ and $1<l_1<\ldots<l_t\leq l$ and
      $s+t\leq l-1$.
\item If $v\in\F^{2l}$ satisfies $v\hat{H}\T=0$ and $v_1\not=0$, then $\wt(v)\geq l+1$.
\end{alphalist}
\end{theo}

\begin{proof}
(a)~$\Leftrightarrow$~(b) is obvious since in case of properness
the Toeplitz structure implies
\[
   T^{i_1,\ldots,i_r}_{j_1,\ldots,j_r}
    =T^{i_1-j_1+1,\ldots,i_r-j_1+1}_{j_1-j_1+1,\ldots,j_r-j_1+1}.
\]
(b)~$\Rightarrow$~(c):
Let $\hat{h}:=T_1+\sum_{j=1}^s \beta_j T_{m_j}$ and
assume to the contrary $\wt(\hat{h})<l-s$.
The assumption implies that $\hat{h}$ consists of at least $s+1$
zero entries, say at the positions $i_1,\ldots,i_{s+1}$.
Then
\begin{equation}\label{e-Abeta}
   T^{i_1,\ldots,i_{s+1}}_{1,m_1,\ldots,m_s}
   \begin{pmatrix}1\\ \beta_1\\ \vdots\\ \beta_s\end{pmatrix}
   =\begin{pmatrix}0\\ 0\\ \vdots\\ 0\end{pmatrix}.
\end{equation}
The superregularity yields $m_{\nu}>i_{\nu+1}$ for some
$\nu\in\{1,\ldots,s\}$, which we can choose to be minimal with this
property.
Then the submatrix $T^{i_1,\ldots,i_{\nu+1}}_{m_{\nu},\ldots,m_s}$ is
identically zero and therefore we obtain from\eqr{e-Abeta} the identity
$T^{i_1,\ldots,i_{\nu}}_{1,m_1,\ldots,m_{\nu-1}}(1,\beta_1,\ldots,\beta_{\nu-1})\T=0$,
a contradiction to superregularity since by minimality of $\nu$
this coefficient matrix is nonsingular.\\
(c)~$\Rightarrow$~(b):
Assume to the contrary that
$\det T^{i_1,\ldots,i_{s+1}}_{1,m_1,\ldots,m_s}=0$ for some indices satisfying
$m_{\nu}\leq i_{\nu+1}$ for $\nu=1,\ldots,s$.
We can assume $s$ to be minimal with this property.
Then there exists $(\beta_0,\beta_1,\ldots,\beta_s)\in\F^{s+1}\backslash\{0\}$ such that
$ T^{i_1,\ldots,i_{s+1}}_{1,m_1,\ldots,m_s}(\beta_0,\ldots,\beta_s)\T=0$.
Minimality of $s$ and the equivalence of~(a) and~(b) imply
$\beta_0\not=0$.
Hence we can take $\beta_0=1$ and\eqr{e-Abeta} is satisfied.
Thus $\wt(T_1+\sum_{j=1}^s\beta_j T_{m_j})\leq l-(s+1)$,
a contradiction.
\\
The properties~(d) and~(e) are simply reformulations of~(c).
\\
The equivalence~(d)~$\Leftrightarrow$~(f) is clear from the structure of
$\hat{H}$ (a linear combination of $T_1$ by the other columns of $\hat{H}$
has to involve the column~$e_1$ and vice versa).
\\
The property~(g) is a reformulation of~(f).
\end{proof}

The equivalence of~(e) and~(g) immediately implies
\begin{cor}\label{C-superreginverse}
If $T\in\F^{l\times l}$ is a superregular lower triangular Toeplitz
matrix, then so is~$T^{-1}$.
\end{cor}

The following lemma is the main step for establishing the existence of
superregular matrices of Toeplitz-structure.

\begin{lemma}\label{L-genericexist}
Let $\F$ be a field and $X_1,\ldots,X_l$ be independent indeterminates
over $\F$. Define the matrix
\[
   A:=\begin{bmatrix}X_1    & 0      & \cdots & 0 \\
                    X_2    & X_1    & \ddots & \vdots  \\
                    \vdots & \ddots &\ddots & 0\\
                    X_l    & \cdots &   X_2  & X_1\end{bmatrix}
   \in\F(X_1,\ldots,X_l)^{l\times l}.
\]
Then $A$ is superregular.
\end{lemma}

\begin{proof}
We proceed by contradiction.
Assume there exists a singular proper submatrix
\[
   \hat{A}:=A^{i_1,\ldots,i_r}_{j_1,\ldots,j_r}.
\]
We can take the size $r$ to be minimal.
Then certainly $r>1$.
By properness we know that $j_{\nu}\leq i_{\nu}$ for
$\nu=1,\ldots,r$.
\\
Notice that for $\mu\leq\nu$ the entry of $A$ at the position $(\nu,\mu)$
is given by $A^{\nu}_{\mu}=X_{\nu-\mu+1}$.
Hence the indeterminate with the largest index appearing in
$\hat{A}$ is $X_{i_r-j_1+1}$. It appears only once in the matrix and that
is in the lower left corner. Thus its coefficient in $\det\hat{A}$ is
$\pm\det\tilde{A}$, where
\[
   \tilde{A}:=A^{i_1,\ldots,i_{r-1}}_{j_2,\ldots,j_r}.
\]
Singularity of $A$ now implies $\det\tilde{A}=0$.
By minimality of $r$ this yields that $\tilde{A}$ is an improper
submatrix of $A$, i.~e.\ there exists an index $\tau\in\{2,\ldots,r\}$
such that $j_{\tau}>i_{\tau-1}$.
Picking $\tau$ minimal we get $i_1<\ldots<i_{\tau-1}<j_{\tau}<\ldots<j_r$
and therefore the first $\tau-1$ rows of $\hat{A}$
have the form
\[
   \begin{bmatrix}*&\cdots&*     &0     &\cdots&0\\
             \vdots&      &\vdots&\vdots&      &\vdots\\
                  *&\cdots& *    &0     &\cdots&0\end{bmatrix},
\]
where the block of possibly nonzero elements consists of $\tau-1$ columns.
Hence $\hat{A}$ is a blocktriangular matrix and we have
\[
  0= \det\hat{A}=\det A^{i_1,\ldots,i_{\tau-1}}_{j_1,\ldots,j_{\tau-1}}
               \det A^{i_{\tau},\ldots,i_r}_{j_{\tau},\ldots,j_r}.
\]
Since both factors are proper minors we get a contradiction to
the minimality of the size~$r$.
\end{proof}

The following consequence is standard.
\begin{theo}\label{T-exist}
For every $l\in\N$ and every prime number~$p$ there exists a finite
field $\F$ of characteristic~$p$ and a superregular
matrix $T\in\F^{l\times l}$ having Toeplitz structure.
\end{theo}

\begin{proof}
Consider the prime field $\F_p$ and the matrix of the previous lemma with entries in
$\F_p(X_1,\ldots,X_l)$.
All its proper
minors are nonzero polynomials in $\F_p[X_1,\ldots,X_l]$.
Over an algebraic closure $\bar{\F}_p$ a point $a:=(a_1,\ldots,a_l)\in\bar{\F}_p^l$
can be found such that none of the minors vanishes at $a$.
Hence the Toeplitz matrix $T$ having $(a_1,\ldots,a_l)\T$ as its first column
is superregular. Since each $a_i$ is algebraic over $\F_p$, the
matrix~$T$ has its entries in a finite field extension $\F$ of $\F_p$.
\end{proof}

In particular, for every size $l\in\N$ there exist superregular Toeplitz matrices
over a field of characteristic~$2$.
Unfortunately, the theorem above is nonconstructive and it is not at all clear what
the minimum field of characteristic~$2$ is to allow a superregular Toeplitz
matrix of given size~$l\times l$.
We present some examples.

\begin{exa}\label{E-superreg2}
\begin{arabiclist}
\item
Using a computer algebra program one checks that the following matrices are
superregular. The first examples are all over prime fields
$\F_p$.

$$
\begin{bmatrix}
1& 0\\
1& 1
\end{bmatrix}\in\F_{2}^{2\times2},
\quad
\begin{bmatrix}
1& 0& 0\\
1& 1& 0\\
2& 1& 1
\end{bmatrix}\in\F_{3}^{3\times3},
\quad
\begin{bmatrix}
1& 0& 0& 0\\
1& 1& 0& 0\\
2& 1& 1& 0\\
1& 2& 1& 1
\end{bmatrix}\in\F_{5}^{4\times 4},
\quad
\begin{bmatrix}
1& 0& 0& 0& 0\\
2& 1& 0& 0& 0\\
1& 2& 1& 0& 0\\
6& 1& 2& 1& 0\\
4& 6& 1& 2& 1
\end{bmatrix}\in\F_{7}^{5\times 5},
$$

$$
\begin{bmatrix}
1& 0& 0& 0& 0& 0\\
2& 1& 0& 0& 0& 0\\
1& 2& 1& 0& 0& 0\\
1& 1& 2& 1& 0& 0\\
3& 1& 1& 2& 1& 0\\
4& 3& 1& 1& 2& 1
\end{bmatrix}\in\F_{11}^{6\times 6},
\quad
\begin{bmatrix}
 1& 0& 0& 0& 0& 0& 0\\
 7& 1& 0& 0& 0& 0& 0\\
13& 7& 1& 0& 0& 0& 0\\
 2&13& 7& 1& 0& 0& 0\\
 1& 2&13& 7& 1& 0& 0\\
 4& 1& 2&13& 7& 1& 0\\
14& 4& 1& 2&13& 7& 1
\end{bmatrix}\in\F_{17}^{7\times 7}.
$$

The following examples represent superregular matrices  over finite fields of
characteristic 2. For this assume that  $\alpha,\beta$ and $\gamma$ satisfy
\[
  \alpha^2+\alpha+1=0,\ \beta^3+\beta+1=0,\text{ and } \gamma^4+\gamma+1=0.
\]
Then the following  matrices represent superregular matrices over
$\F_4$, $\F_8$ and $\F_{16}$ respectively.
\[
     \begin{bmatrix}1& & \\ \alpha&1& \\
             1&\alpha&1\end{bmatrix}\in\F_{2^2}^{3\times3},\quad
     \begin{bmatrix}1& & & & \\ \beta&1& & & \\
                       \beta^3&\beta&1& & \\ \beta&\beta^3&\beta&1& \\
                       1&\beta&\beta^3&\beta&1
     \end{bmatrix}\in\F_{2^3}^{5\times5},\quad
     \begin{bmatrix}1& & & & & \\ \gamma&1& & & & \\ \gamma^5&\gamma&1& & & \\
                     \gamma^5&\gamma^5&\gamma&1& & \\
                     \gamma&\gamma^5&\gamma^5&\gamma&1& \\
                     1&\gamma&\gamma^5&\gamma^5&\gamma&1
     \end{bmatrix}\in\F_{2^4}^{6\times6}.
\]
Assume $\epsilon,\omega$ satisfy
\[
   \epsilon^5+\epsilon^2+1=0 \text{ and } \omega^6+\omega+1=0.
\]
Then the following matrices represent superregular matrices over
$\F_{32}$ and $\F_{64}$ respectively.
\[
     \begin{bmatrix}1& & & & & & \\ \epsilon&1& & & & & \\
                      \epsilon^6&\epsilon&1& & & & \\
                      \epsilon^9&\epsilon^6&\epsilon&1& & & \\
                      \epsilon^6&\epsilon^9&\epsilon^6&\epsilon&1& & \\
                      \epsilon&\epsilon^6&\epsilon^9&\epsilon^6&\epsilon&1& \\
                      1&\epsilon&\epsilon^6&\epsilon^9&\epsilon^6&\epsilon&1
     \end{bmatrix}\in\F_{2^5}^{7\times7},\quad
     \begin{bmatrix}1&0&0&0&0&0&0&0\\ \omega&1&0&0&0&0&0&0\\
                   \omega^9&\omega&1&0&0&0&0&0\\\omega^{33}&\omega^9&\omega&1&0&0&0&0\\
                   \omega^{33}&\omega^{33}&\omega^9&\omega&1&0&0&0\\
                   \omega^9&\omega^{33}&\omega^{33}&\omega^9&\omega&1&0&0\\
                   \omega&\omega^9&\omega^{33}&\omega^{33}&\omega^9&\omega&1&0\\
                   1&\omega&\omega^9&\omega^{33}&\omega^{33}&\omega^9&\omega&1
     \end{bmatrix}\in\F_{2^6}^{8\times8}.
\]

Notice that the matrices above have even more symmetry than required.
One can easily show that there is no superregular $4\times4$-matrix over
$\F_4$ of general Toeplitz structure.
However, the above suggests to ask whether one can find for every $l\geq5$ a superregular
$l\times l$-Toeplitz matrix over $\F_{2^{l-2}}$.
\item
In the appendix we prove that for every $n\in\N$ the proper minors of the Toeplitz-matrix
\[
     T_n:=\begin{bmatrix}
        \genfrac{(}{)}{0pt}{}{n-1}{0}&           &  &  \\
       \genfrac{(}{)}{0pt}{}{n-1}{1}&\genfrac{(}{)}{0pt}{}{n-1}{0}&        &  \\
              \vdots & \ddots      &\ddots  &\rule[-.3cm]{0cm}{1cm}  \\
       \genfrac{(}{)}{0pt}{}{n-1}{n-1}& \cdots &\genfrac{(}{)}{0pt}{}{n-1}{1}&\genfrac{(}{)}{0pt}{}{n-1}{0}
        \end{bmatrix}\in\Z^{n\times n}
\]
are all positive.
Hence for each $n\in\N$ there exists a smallest prime number $p_n$ such that $T_n$
is superregular over the prime field $\F_{p_n}$.
One can check that
\[
   p_2=2,\ p_3=5,\ p_4=7,\ p_5=11,\ p_6=23,\ p_7=43.
\]
\end{arabiclist}
\end{exa}

Now we can establish the existence of strongly MDS codes in the following sense.

\begin{theo}\label{T-str.MDSn}
For every $n,\,\delta\in\N$ and every prime number~$p$ there exists a
strongly MDS code with parameters $(n,n-1,\delta)$ over a suitably large field of
characteristic~$p$.
\end{theo}

The proof of this theorem is rather long and technical and
because of this reason it is put into the appendix.

There is of course the natural question if strongly MDS
convolutional codes and codes with maximum distance profile exist
for all parameters $(n,k,\delta)$. We strongly believe so. The
section showed that such codes exist for all parameters
$(n,k,\delta)$ with $k=n-1$. For all small values of
$(n,k,\delta)$ we have found strongly MDS
convolutional codes and codes with maximum distance profile
making computer searches. In the next section we present a series
of examples of such codes found through computer searches. Based
on this wealth of data we conjecture:
\begin{conjecture}
  For all $n>k>0$ and for all $\delta\geq 0$ there exists an
  $(n,k,\delta)$ code over a sufficiently large field which is
  both strongly MDS and has a maximum distance profile.
\end{conjecture}

\section{Examples}
\setcounter{equation}{0}                      \label{Sec4}

In this section we will present some examples of strongly MDS codes with small
parameters.
The first set of examples is constructed according to the proof of Theorem~\ref{T-str.MDSn}
by utilizing the superregular matrices in Example~\ref{E-superreg2}.

\begin{exa}\label{E-MDS}
Recall the first part of the proof of Theorem~\ref{T-str.MDSn}.
\begin{arabiclist}
\item We can construct strongly MDS $(2,1,\delta)$-codes once a
  $\tau\times\tau$ superregular matrix, where $\tau=2\delta+1$,
  is available. Thus, the $5\times 5$ and $7\times 7$ matrices
  given in Example~\ref{E-superreg2}(1) lead to the strongly MDS
  $(2,1,2)$-code over $\F_{8}$ (where $\beta^3+\beta+1=0$) with
  parity check matrix
      \[
          H=[a,\,b]=[1+\beta^2D+\beta^5D^2,\,1+\beta^4D+\beta^5D^2]\in\F_{8}[D]^{2}
      \]
      and to the strongly MDS $(2,1,3)$-code over $\F_{32}$ (where $\epsilon^5+\epsilon^2+1=0$)
      with parity check matrix
      \[
         H=[a,\,b]=[1+\epsilon^{18}D+\epsilon^{11}D^2+\epsilon^{29}D^3,\,
                  1+D+\epsilon^{27}D^2+\epsilon^{18}D^3]\in\F_{32}[D]^{2}.
      \]
      Indeed, one checks that
      \[
        \frac{1+\beta^4D+\beta^5D^2}{1+\beta^2D+\beta^5D^2}=1+\beta D+\beta^3D^2+\beta
        D^3+D^4+\text{ higher powers}
      \]
      and
      \[
        \frac{1+D+\epsilon^{27}D^2+\epsilon^{18}D^3}{1+\epsilon^{18}D+
        \epsilon^{11}D^2+\epsilon^{29}D^3}
        =\!1+\!\epsilon D\!+\!\epsilon^6 D^2\!+\!\epsilon^9 D^3\!+\!\epsilon^6 D^4\!+\!\epsilon
        D^5\!+\!D^6+\!\!\text{ higher powers.}
      \]
      Hence the free distance of the two codes above is~$6$ (resp.~$8$), and this is
      also the $4$th (resp.~$6$th) column distance.
\item Using the $8\times8$-superregular matrix of
      Example~\ref{E-superreg2}(1), one can construct a strongly MDS
      $(3,2,2)$-code over $\F_{64}$.
      Hence the code has free distance equal to its $3$rd column distance, and this value
      is~$5$. Using the construction of the proof of
      Theorem~\ref{T-str.MDSn} and going through some tedious calculations in the field $\F_{64}$
      (where $\omega^6+\omega+1=0$) one finally arrives at the parity check matrix
      \[
         H=[1+\omega^{57}D+\omega^{62}D^2,\,\omega+\omega^{44}D+\omega^{54}D^2,\,
            1+\omega^{17}D+\omega^{21}D^2]\in\F_{64}^3.
      \]
\item A strongly MDS $(4,3,1)$-code has free distance~$3$ and this is identical with
      the first column distance. It can be obtained from a $6\times6$-superregular
      matrix using the construction of the proof of Theorem~\ref{T-str.MDSn}. Indeed,
      the matrix
      \[
        \hat{H}=\left[\!\!\begin{array}{cc|cccccc} 1&0&\gamma^5&\gamma&1&0&0&0\\
        0&1&1&\gamma&\gamma^5&\gamma^5&\gamma&1\end{array}\!\!\right]
      \]
      has been obtained from the superregular Toeplitz matrix of Example~\ref{E-superreg2}(1) and
      thus it satisfies property~(b) of Theorem~\ref{T-strongMDS}.
      Hence a parity check matrix of a strongly MDS
      $(4,3,1)$-code over $\F_{16}$ (where $\gamma^4+\gamma+1=0$)
      is given by
      \[
          H=[1,\gamma^5+D,\gamma+\gamma D,1+\gamma^5D]\in\F_{16}[D]^{4}.
      \]
\item Of course, not every MDS code is strongly MDS.
      For instance, the code with parity check matrix
      $H=[10+3D+2D^2,4+2D+D^2]\in\F_{11}[D]^2$ is an MDS code,
      but not strongly MDS.
      In this example, the MDS property follows from the fact, that this code is the result of
      the construction of MDS codes as presented in~\cite{sm01a}.
      However, a $(2,1,1)$-code is strongly MDS iff it is an MDS code.
      This can be checked directly by using Theorem~\ref{T-strongMDS} and
      the fact that for the (basic) parity check matrix $[a_0+a_1D,b_0+b_1D]$ of an
      MDS code all coefficients as well as $a_0b_1-a_1b_0$ are nonzero.
\end{arabiclist}
\end{exa}

The next series of examples has been found by completely
different methods.  They are all cyclic convolutional codes in
the sense of~\cite{gl02u,gl02p,pi76,ro79}.  In those papers
convolutional codes having some additional algebraic structure
are being investigated. This additional structure is a
generalization of cyclicity of block codes but is a far more
complex notion for convolutional codes. In particular cyclicity
of convolutional codes does {\em not\/} mean invariance under the
cyclic shift in $\FDlaurent^n$.  We will not go into the details but
rather refer to~\cite{gl02u,gl02p}.  However, in order to
understand and test the following examples there is no need in
understanding the concept of cyclicity for convolutional codes
since below we provide all information needed to specify the
codes.  We present the generator matrices and also provide all
column distances; they have been computed with a computer algebra
program.  All matrices given below are minimal basic in the sense
of Definition~\ref{D-basic.minimal}.  We would like to mention
that just like for cyclic block codes, the length of the code and
the characteristic of the field have to be coprime. Therefore,
only codes with odd length are given below.

One should note that most of the following codes exist over comparatively smaller alphabets
than the examples of~\ref{E-MDS}.
However, we don't know any general construction for strongly MDS cyclic convolutional
codes yet.
But the abundance of (small) examples suggests that such a construction
might be possible and might lead to smaller alphabets for given parameters than the
construction of the last section.
We will leave this as an open question for future research.

\begin{exa}\label{E-strMDSCCC}
\begin{arabiclist}
\item A strongly MDS $(3,1,1)$-code over $\F_4$:
      \[
        G=[\alpha+\alpha D,\,\alpha^2+\alpha D,\,1+\alpha D].
      \]
      The column distances are $d^c_0=3,\,d^c_1=5,\,d^c_j=6$ for $j\geq2$.
\item A strongly MDS $(3,1,2)$-code over $\F_{16}$ (where $\beta^4+\beta+1=0$):
      \[
        G=[\beta+\beta D+D^2,\,\beta^6+\beta D+\beta^{10}D^2,\,\beta^{11}+\beta D+\beta^5 D^2].
      \]
      The column distances are $d^c_0=3,\,d^c_1=5,\,d^c_2=7,\,d^c_j=9$ for $j\geq3$.
\item A strongly MDS $(3,2,2)$-code over $\F_{16}$:
      \[
        G=\begin{bmatrix}
            \beta^5+\beta^4 D&\beta^3+\beta^8D&\beta^9+\beta^2D\\
            \beta^9+\beta^{12}D&\beta^5+\beta^{14}D&\beta^3+\beta^3D
          \end{bmatrix}.
      \]
      The column distances are $d^c_0=2,\,d^c_1=3,\,d^c_2=4,\,d^c_j=5$ for $j\geq3$.
\item A strongly MDS $(5,1,1)$-code over $\F_{16}$:
      \[
        G=[\beta+\beta D,\,\beta^{13}+\beta^{10}D,\,\beta^{10}+\beta^4D,\,\beta^7+\beta^{13}D,\,
           \beta^4+\beta^7D].
      \]
      The column distances are $d^c_0=5,\,d^c_1=9,\,d^c_j=10$ for $j\geq2$.
\item A strongly MDS $(5,1,2)$-code over $\F_{16}$:
      \begin{align*}
        &G=[\beta+\beta^4 D+\beta D^2,\,\beta^{7}+\beta D+\beta^{10}D^2,\,
            \beta^{13}+\beta^{13}D+\beta^4 D^2,\\
        &\mbox{}\hspace*{17em}
            \beta^4+\beta^{10}D+\beta^{13}D^2,\,\beta^{10}+\beta^7D+\beta^7D^2].
      \end{align*}
      The column distances are $d^c_0=5,\,d^c_1=9,\,d^c_2=13,\,d^c_j=15$ for $j\geq3$.
\item A strongly MDS $(5,2,2)$-code over $\F_{16}$:
      \[
        G=\begin{bmatrix}
           \beta+\beta
           D&\beta^{13}+\beta^{10}D&\beta^{10}+\beta^4D&\beta^7+\beta^{13}D&
           \beta^4+\beta^7D\\
           1+\beta^5D&\beta^3+\beta^{11}D&\beta^6+\beta^2D&\beta^9+\beta^8D&\beta^{12}+\beta^{14}D
           \end{bmatrix}.
      \]
      The column distances are $d^c_0=4,\,d^c_1=7,\,d^c_j=9$ for $j\geq2$.
\item A strongly MDS $(7,1,1)$-code over $\F_{8}$ (where $\gamma^3+\gamma+1=0$):
      \[
        G=[\gamma+\gamma D,\,\gamma^3+D,\,\gamma^5+\gamma^6D,\,1+\gamma^5D,\,
           \gamma^2+\gamma^4D,\,\gamma^4+\gamma^3D,\,\gamma^6+\gamma^2D].
      \]
      The column distances are $d^c_0=7,\,d^c_1=13,\,d^c_j=14$ for $j\geq2$.
\item A strongly MDS $(7,1,2)$-code over $\F_{8}$:
      \begin{align*}
       & G=[\gamma^2+\gamma D+D^2,\,\gamma^5+\gamma^3D+\gamma^6D^2,\,\gamma+\gamma^5D+\gamma^5D^2,\,
           \gamma^4+D+\gamma^4D^2,\\
       &\mbox{}\hspace*{13em}1+\gamma^2D+\gamma^3D^2,\,\gamma^3+\gamma^4D+\gamma^2D^2,\,
           \gamma^6+\gamma^6D+\gamma D^2].
      \end{align*}
      The column distances are $d^c_0=7,\,d^c_1=13,\,d^c_2=18,\,d^c_j=21$ for $j\geq3$.
\item It is worth being mentioned that there does not exist even an MDS $(7,2,2)$-code
      over $\F_8$, since the generalized Singleton bound in this case is~$13$, but due to the Griesmer
      bound (see~\cite[p.~133]{jo99} for the binary case) the parameters of an $(n,k,\delta)$-code
      over $\F_q$ with memory~$m$ and distance~$d$ satisfy
      \[
        \sum_{l=0}^{k(m+i)-\delta-1}\Big\lceil \frac{d}{q^l}\Big\rceil\leq n(m+i)
        \text{ for all } i\in\N_0.
      \]
      Hence a $(7,2,2)$-code over $\F_8$ with memory~$1$ has at most distance~$12$.
      The inequality applied to $i=1$ shows that
      the field size has to be at least~$13$ in order to allow
      the existence of an MDS $(7,2,2)$-code.
\end{arabiclist}
One should notice that the codes in Example~\ref{E-strMDSCCC}(1)~--~(7) are not only strongly MDS but
also have {\em all\/} column distances being optimal in the sense that they reach the
upper bound given in Proposition~\ref{P-dcj.bound}. In particular
they also have a maximum distance profile in the sense of
Definition~\ref{max-Profile}.
For the $(7,1,2)$-code in~(8), only the second column distance is not optimal, but rather
one less than the upper bound, which is~$19$ in this case.
\end{exa}

\section{The Dual of a Strongly MDS Code}
\setcounter{equation}{0}                                   \label{Sec5}
In this section we will present some results concerning the dual code of a strongly
MDS code. The main result shows that a convolutional code has a
maximum distance profile if and only if its dual has this
property. This then implies  for certain parameters that a code
is strongly MDS if and only if its dual has this property. These
results are very appealing as it generalizes the situation for
block codes.

Recall that if
\[
  \C=\{uG\mid u\in\FDlaurent^k\}=\{v\in\FDlaurent^n\mid v H\T=0\}\subseteq\FDlaurent^n
\]
is an $(n,k,\delta)$-code with generator matrix
$G\in\F[D]^{k\times n}$ and parity check matrix
$H\in\F[D]^{(n-k)\times n}$, then the dual code, defined as
\[
  \C^{\perp}=\{w\in\FDlaurent^n\mid wv\T=0\text{ for all }v\in\C\},
\]
is given by
\[
  \C^{\perp}=\{uH\mid u\in\FDlaurent^{n-k}\}=\{w\in\FDlaurent^n\mid wG\T=0\}
\]
and thus an $(n,n-k,\delta)$-code.  In contrast to the block code
situation almost nothing is known about the relation between the
distances of a code and its dual.  In particular, it has been
shown in~\cite{sh77a} that no MacWilliams identity relating the
weight distributions of~$\C$ and~$\C^{\perp}$ exists.  In block
code theory a very simple relation between the distances of a
code and its dual is given in the case of MDS codes. In fact,
if~$\C$ is an MDS $(n,k)$-block-code, then the dual~$\C^{\perp}$
is an MDS $(n,n-k)$-code, see~\cite[Ch.~11, \S2]{ma77} and very
specific knowledge on the weight enumerator and its dual is
known~\cite[Ch.~11]{ma77}.  Therefore, it is quite natural to
investigate whether the dual of an MDS (or strongly MDS)
convolutional code is MDS (or strongly MDS), too.  Unfortunately,
this is in general not the case.

\begin{exa}\label{E-Cdual1}
In general the dual of a strongly MDS code is not even an MDS code.
This can be seen from the dual of the code given in Example~\ref{E-MDS}(3).
The dual has generator matrix
$G=[1,\gamma^5+D,\gamma+\gamma D,1+\gamma^5D]\in\F_{16}[D]^{4}$ which obviously has
weight less than the generalized Singleton bound~$8$ (see Theorem~\ref{T-MDS}).
\end{exa}

As we will show next the property of maximum distance profile
carries over under dualization. In addition,
 for specific code parameters the strong MDS property carries over to
the dual code as well. To this end, recall from Definition~\ref{D-strongly.MDS} that an
$(n,k,\delta)$-code is strongly MDS if the $M$th column distance attains
the generalized Singleton
bound where $M=\lfloor\frac{\delta}{k}\rfloor+\lceil\frac{\delta}{n-k}\rceil$.
Thus the dual code~$\C^{\perp}$ is MDS if the $\hat{M}$th column distance attains the
generalized Singleton bound where
$\hat{M}=\lfloor\frac{\delta}{n-k}\rfloor+\lceil\frac{\delta}{k}\rceil$.
Obviously, these two numbers differ by one when $k$ divides $\delta$
but~$n-k$ does not or vice versa. What remains equal for both the code and its
dual is the quantity
$L=\Big\lfloor\frac{\delta}{k}\Big\rfloor+\Big\lfloor\frac{\delta}{n-k}\Big\rfloor$
used in Definition~\ref{max-Profile} where we introduced the
concept of maximum distance profile.

Before we state the main results we need a technical lemma.
\begin{lemma}                         \label{L-minors}
Let $A\in\F^{k\times n}$ and $B\in\F^{n\times(n-k)}$ such that
\[
 AB=0\text{ and }\rank A=k,\ \rank B=n-k.
\]
Then the following are equivalent:
\begin{alphalist}
\item the $k\times k$-submatrix of~$A$ consisting of the columns with indices
      $1\leq t_1<\ldots<t_k\leq n$ is singular,
\item The $(n-k)\times(n-k)$-submatrix of~$B$ obtained by taking the rows with
      indices in $\{1,\ldots,n\}\backslash\{t_1,\ldots,t_k\}$ is singular.
\end{alphalist}
\end{lemma}

\begin{proof}
Without loss of generality assume $(t_1,\ldots,t_k)=(1,\ldots,k)$
and partition $A=(A_1\ A_2)$, where $A_1$ is the $k\times k$
submatrix under consideration. If $A_1$ is invertible then
$$
\ker A=\mathrm{colspan}_\F\left( A_1^{-1}A_2 \atop-I_{n-k} \right)
=\mathrm{colspan}_\F(B).
$$
This shows that the bottom $(n-k)\times(n-k)$-submatrix of~$B$ is
invertible.
\end{proof}

This lemma, in conjunction with Theorem~\ref{Alg-criterion}
immediately gives an algebraic criterion for maximum distance
profile codes in terms of a parity check matrix.

\begin{theo}                        \label{Alg-criterion2}
 Let $H=\sum_{j=0}^\mu H_jD^j$ be the parity check  matrix of an
 $(n,k,\delta)$-code.
  Let $L$ be defined as in\eqr{e-L} and let
 \begin{equation}
\begin{array}{rcl}
 H^c_L: &=& \begin{bmatrix}
                 H_0&    &   &      \\
                 H_1& H_0&   &      \\
                \vdots&\vdots&\ddots&\\
                 H_L&H_{L-1} &\ldots& H_0
              \end{bmatrix}\in\F^{(L+1)(n-k)\times(L+1)n}.
\end{array}
\end{equation}
Then $H$ represents a maximum distance profile code if and only
if every $(L+1)(n-k)\times(L+1)(n-k)$ full-size minor formed from
the columns with indices $1\leq i_1< \ldots < i_{(L+1)(n-k)}$,
where $i_{s(n-k)}\leq sn$ for $s=1,\ldots,L$, is nonzero.
\end{theo}

\begin{proof}
Let the code have generator matrix~$G$ as given in\eqr{e-G}.
Recall that $G^c_L(H^c_L)\T=0$ and both factors have full rank.
By Theorem~\ref{Alg-criterion} the code has maximum distance profile
if and only if every full size minor $G^c_L$
formed from the columns $1\leq j_1< \ldots < j_{(L+1)k}$, where
$j_{sk+1}>sn$ for $s=1,\ldots,L$, is nonzero.
Now the complimentary minors of $H^c_L$ have indices
$1\leq i_1< \ldots < i_{(L+1)(n-k)}$ satisfying
$i_{s(n-k)}\leq sn$ for $s=1,\ldots,L$.
Thus Lemma~\ref{L-minors} completes the proof.
\end{proof}

With this we have a nice duality result:

\begin{theo}                         \label{T-strMDSdual}
An $(n,k,\delta)$-code $\C\subseteq\FDlaurent^n$ has a maximum
distance profile if and only if the dual code
$\C^{\perp}\subseteq\FDlaurent^n$ has this property.
\end{theo}

\begin{proof}
Let $\C$ have generator matrix $G$ and parity check matrix $H$
as given in\eqr{e-G} and\eqr{e-H}.  Assume $\C$ has a maximum
distance profile. By Theorem~\ref{Alg-criterion2} every
$(L+1)(n-k)\times(L+1)(n-k)$ full-size minor formed from the
columns of $H^c_L$ with indices $1\leq i_1< \ldots <
i_{(L+1)(n-k)}$, where $i_{s(n-k)}\leq sn$ for $s=1,\ldots,L$,
is nonzero.

Consider now the dual code $\C^{\perp}$ which
is defined as the rowspace of the $(n-k)\times n$ matrix
$H$. It follows from\eqr{e-dcj} that the $L$th column distance
of the dual code $\C^{\perp}$ is given by
$$
 \hat{d}^c_L=\min\big\{\wt\big((u_L,\ldots,u_0)H^c_L\big)\,\big|\,
                u_i\in\F^{n-k},\,u_0\not=0\big\}.
$$
Taking the reversed ordering into account we obtain from
Theorem~\ref{Alg-criterion} that the dual
code $\C^{\perp}$ has maximum distance profile as well.
\end{proof}

\begin{cor}                             \label{dual-special}
  When both~$k$ and~$n-k$ divide~$\delta$ then an  $(n,k,\delta)$-code
  $\C\subseteq\FDlaurent^n$ is strongly MDS if and only if
  $\C^{\perp}\subseteq\FDlaurent^n$ has this property.
\end{cor}

\begin{proof}
{}From $k\mid\delta$ and $(n-k)\mid\delta$ it follows that $L=M$
and $d^c_M= (n-k)\big(\frac{\delta}{k}+1\big)+\delta+1$, the
generalized Singleton bound of the code  $\C$ and  $\hat{d}^c_M=
k\big(\frac{\delta}{n-k}+1\big)+\delta+1$, the generalized
Singleton bound of the dual code  $\C^{\perp}$.
\end{proof}

The result above gives us another class of strongly MDS codes by
dualizing Theorem~\ref{T-str.MDSn}.

\begin{cor}\label{C-strMDSdual}
For every $n,\,\delta\in\N_0$ such that $(n-1)\mid\delta$ and every prime number~$p$
there exists a strongly MDS $(n,1,\delta)$-code over some suitably large field
of characteristic~$p$.
\end{cor}
\begin{exa}\label{E-dualMDS}
\begin{alphalist}
\item Corollary~\ref{dual-special} tells us that the duals of the $(2,1,\delta)$-codes
      given in Example~\ref{E-MDS}(1) are strongly MDS. But this is obviously so,
      since they are --- up to ordering --- identical to the given codes.
\item Dualizing the code of Example~\ref{E-MDS}(2) gives us a strongly MDS
      $(3,1,2)$-code with generator matrix
      \[
          G=[1+\omega^{57}D+\omega^{62}D^2,\,\omega+\omega^{44}D+\omega^{54}D^2,\,
            1+\omega^{17}D+\omega^{21}D^2]\in\F_{64}^3.
      \]
\item Dualizing the codes given in Example~\ref{E-strMDSCCC}(2) and~(3) we obtain
      another two strongly MDS codes with generator matrices
      \[
         H_1=\begin{bmatrix}
             1 & \beta D+\beta^9 & \beta^6D+\beta^8\\
             \beta^{14}D & \beta^7D+\beta^6 &
             \beta^8D+\beta \end{bmatrix}\in\F_{16}^{2\times3}
      \]
      and
      \[
        H_2=[D^2+D+\beta^2, \beta^{10}D^2+D+\beta^7, \beta^5D^2+D+\beta^{12}]\in\F_{16}^3.
      \]
      It is known that these codes are also cyclic convolutional codes in the sense
      of~\cite{gl02u}, see~\cite[Thm.~7.5]{gl02u}.
\end{alphalist}
\end{exa}

Finally we would like to mention that even in the case where $k\mid\delta$ and $(n-k)\mid\delta$,
the dual of an MDS code is not MDS in general.
An example is given by the following code.

\begin{exa}\label{E-MDSdual}
The $(3,1,2)$-code $\C\subseteq\FDlaurent^3$, where $\F=\F_{16}$,  with generator matrix
\[
   G=[1+\beta D+\beta^4D^2,\ \beta^{10}+\beta^2D+\beta^4D^2,\ \beta^8+\beta^5D+D^2]
\]
and parity check matrix
\[
   H=\begin{bmatrix}1 & \beta^{14}D+\beta^2 & \beta^3D+\beta^3\\
        \beta D & \beta^{11}D+\beta^8 & \beta^{10}D+\beta^{10}
     \end{bmatrix}
\]
is an MDS code, but not strongly MDS. It satisfies $d^c_3=8$ and $d^c_4=9$.
The dual code generated by~$H$ is not MDS. Its distance is~$4$.
\end{exa}

\section{Decoding Strongly MDS Codes}
\setcounter{equation}{0}                           \label{Sec6}
The codes discussed in the previous section have the property that they
allow a very good feedback decoding~\cite{Ro68a} if the error distribution is
reasonably mild.

Let us briefly recall the concept of feedback decoding.
Assume the codeword $v=\sum_{j\geq0}v_j D^j\in\C$ has been sent
and the word $\hat{v}=\sum_{j\geq0}\hat{v}_jD^j\in\FDlaurent^n$ has been
received.
Write $\hat{v}=v+\epsilon$, where
$\epsilon=\sum_{j\geq0}\epsilon_jD^j$ is the error vector.
In the $j$-th cycle of feedback decoding one corrects $\hat{v}_j$
(hence estimates $\epsilon_j$) and then feeds back
this information into the decoding
algorithm before proceeding with the next decoding step.
It should be intuitively clear that the next step will benefit from the resetting
$\hat{v}\leftarrow \hat{v}-\epsilon_jD^j$.
As for the decoding step itself one estimates $\epsilon_j$ upon the knowledge of
the received sequence $\hat{v}_j,\ldots,\hat{v}_{j+l}$.
The length $l+1$, of course, depends on the distance properties of the code.
This estimate will be correct if not too many errors
have occurred on this string.

In the sequel we will show that strongly MDS codes of rate $\frac{n-1}{n}$ have very good
error correcting capabilities in terms of the maximum number of errors
acceptable on each string without jeopardizing correct decoding.
The basis of the feedback decoding algorithm is the following simple reformulation
of the distance properties for the parity check matrices.

\begin{prop}\label{P-HcM.error}
Let $\C\subseteq\FDlaurent^n$ be a strongly MDS $(n,n-1,\delta)$-code and put
$M:=\floor{\frac{\delta}{n-1}}+\delta$.
Let $H^c_M\in\F^{(M+1)\times(M+1)n}$ be the $M$-th parity check matrix of
$\C$ and
$\epsilon:=(\epsilon_0,\ldots,\epsilon_M),\,
 \hat{\epsilon}:=(\hat{\epsilon}_0,\ldots,\hat{\epsilon}_M)\in\F^{M+1}$.
Assume
\[
  \epsilon(H^c_M)\T=\hat{\epsilon}(H^c_M)\T
  \text{ and }
  \wt(\epsilon),\,\wt(\hat{\epsilon})\leq\frac{M+1}{2}.
\]
Then
\begin{alphalist}
\item $\epsilon_0=\hat{\epsilon}_0$,
\item if $M$ is even, then additionally $\epsilon_1=\hat{\epsilon}_1$.
\end{alphalist}
Notice that $M$ is even for codes with rate $1/2$.
\end{prop}

\begin{proof}
Put $\eta:=(\eta_0,\ldots,\eta_M)$ where
$\eta_j=\epsilon_j-\hat{\epsilon}_j$ for all $0\leq j\leq M$.
Then $\eta(H^c_M)\T=0$ and $\wt(\eta)\leq M+1$.
Thus Proposition~\ref{P-dcj} yields $\eta_0=0$.
As for (b) notice that if~$M$ is even, then
$\wt(\epsilon),\,\wt(\hat{\epsilon})\leq\frac{M}{2}$ and therefore
$\wt(\eta)=\wt(\eta_1,\ldots,\eta_M)\leq M$.
Now $\eta_0=0$ implies $(\eta_1,\ldots,\eta_M)(H^c_{M-1})\T=0$ and
Proposition~\ref{P-dcj} together with Corollary~\ref{C-optdistances}
finishes the proof.
\end{proof}

Observe that the proposition above says that the list of $M+1$ consecutive
syndromes determines uniquely the error in the first position. This can be
iterated and leads to the following algorithm, which at least works
reasonably well for small codes. The question how to practically compute the error
in the first position from the syndrome vector for large codes will be
addressed at the end of this section.

We will make use of the notation in Remark~\ref{R-GcjHcj}.

\begin{theo}\label{T-dec.alg}
Let $\C\subseteq\FDlaurent^n$ be a strongly MDS $(n,n-1,\delta)$-code with
parity check matrix $H\in\F[D]^{1\times n}$ and
$H^c_M\in\F^{(M+1)\times(M+1)n}$ as its $M$-th sliding parity check matrix.
Assume the codeword $v\in\C$ has been sent and the word $\hat{v}\in\FDlaurent^n$
has been received.
Without loss of generality assume $\delay{v},\,\delay{\hat{v}}\geq0$.
Put $\hat{v}=v+\epsilon$, thus
$\epsilon\in\FDlaurent^n$ is the error vector and assume that
any sliding window of length $(M+1)n$ contains at most $\frac{M+1}{2}$
errors, i.~e.
\begin{equation}\label{e-slide.error}
  \wt(\epsilon_{[j,j+M]})\leq\frac{M+1}{2}\text{ for all }j\geq0.
\end{equation}
Then the following algorithm will decode $\hat{v}$ correctly, i.~e.\
for each $j=0,1,2,\ldots$ we have $\hat{v}_{[0,j]}=v_{[0,j]}$
after the $j$-th cycle:
\\
Put $j:=-1$.
\begin{algo}
\item Put $j:=j+1$.
\item Compute the syndrome vector
      $S:=(\hat{v}H\T)_{[j,j+M]}$.
\item From the syndrome vector $S$ determine the unique $\eta_0\in\F^n$,
      such that $S=\eta(H^c_M)\T$ for some $\eta=(\eta_0,\ldots,\eta_M)
      \in\F^{(M+1)n}$ satisfying $\wt(\eta)\leq\frac{M+1}{2}$.
\item Put $\hat{v}:=\hat{v}-\eta_0D^j$.
\item Go to Step~1.
\end{algo}
\end{theo}

\begin{rem}\label{R-slide.error}
  For illustration purposes assume $n=2$, i.e. the rate is $1/2$.
  Theorem~\ref{T-dec.alg} then states that a strongly MDS
  $(2,1,\delta)$-code can be correctly decoded as long as there
  are no more than $\delta$ errors in any sliding window of
  length $4\delta +2$. This has to be compared with a MDS block
  code of rate $k/n$ where $n=2k=4\delta+2$ which is capable of
  decoding correctly $\delta$ errors in any slotted window of
  length $n$.  Similar comparisons can be made for different
  values of $n$.
\end{rem}

\begin{proof}
We first have to show the existence of~$\eta$ as required in Step~3 and the uniqueness
of~$\eta_0$.
In order to do so fix some $j\geq0$.
It is easy to see that for all $w\in\FDlaurent^n$ with $\delay{w}\geq0$ one has
\begin{equation}\label{e-syndrome}
  (w H\T)_{[j,j+M]}=w_{[j,j+M]}(H^c_M)\T+w_{[0,j-1]}\H_j\T
\end{equation}
where
\[
   \H_j=\begin{bmatrix}H_j&\cdots&H_1\\H_{j+1}&\cdots&H_2\\ \vdots& &\vdots\\
                       H_{j+M}&\cdots&H_{M+1}\end{bmatrix}.
\]
Since, due to the previous decoding steps we have $\hat{v}_{[0,j-1]}=v_{[0,j-1]}$,
which is the correct codeword sequence, we get
\[
   0=(vH\T)_{[j,j+M]}=\hat{v}_{[j,j+M]}(H^c_M)\T-\epsilon_{[j,j+M]}(H^c_M)\T
      +\hat{v}_{[0,j-1]}(\H_j)\T.
\]
Again with\eqr{e-syndrome} this yields
\[
   S=(\hat{v}H\T)_{[j,j+M]}=\epsilon_{[j,j+M]}(H^c_M)\T
\]
and the assumption\eqr{e-slide.error} together with Proposition~\ref{P-HcM.error} establish the
existence of~$\eta$ as well as the uniqueness of~$\eta_0$ as required in~Step~3.
\\
It follows directly from the above that $\eta_0=\epsilon_j$, where $\eta_0$ is
computed in Step~3 of  the $j$-th cycle. Thus we have $\hat{v}_{[0,j]}=v_{[0,j]}$
after the $j$-th cycle.
\end{proof}

\begin{rem}\label{R-feedb.dec}
One might wonder how the algorithm above is related to the total error
correcting bound $t:=\lfloor\frac{\dfree-1}{2}\rfloor$ of the code. First
notice that $\frac{M+1}{2}=\frac{\dfree-1}{2}=\frac{d^c_M-1}{2}$. From
this it follows that for each received word $\hat{v}$ there exists at most
one codeword $v\in\C$ such that $v-\hat{v}$ satisfies\eqr{e-slide.error}.
This codeword, of course, is then the result of the algorithm above.
However, it might happen that there are two
codewords $v_1,\,v_2\in\C$ such that the total distances satisfy
$\wt(\hat{v}-v_1)=\wt(\hat{v}-v_2)=d(\hat{v},\C):=\min\{\wt(\hat{v}-v)\mid
v\in\C\}$. Hence~$v_1$ and~$v_2$ have equally close  distance
to~$\hat{v}$ when considered over the total length $[0,\infty)$. This of
course can happen only if
$d(\hat{v},\C)>\lfloor\frac{\dfree-1}{2}\rfloor$.
From the above we know that at most one of these codewords can have an error
vector satisfying\eqr{e-slide.error}.
In this situation the decoding algorithm
will try to successively minimize $\wt\big((\hat{v}-v)_{[j,j+M]}\big)$
over all codewords~$v\in\C$ and $j\geq0$.
\\
This situation arises for instance for the strongly MDS code $\C$ with
parity check matrix
\[
    H=[1+\beta^2D+\beta^5D^2,\,1+\beta^4D+\beta^5D^2]\in\F_{2^3}[D]^{2}
    \ (\text{where }\beta^3+\beta+1=0)
\]
given in Example~\ref{E-MDS}(1) and having free distance $\dfree=6$.
In this case the received word
\[
   \hat{v}=(\beta D+\beta^5 D^4,\, \beta^3 D^2+\beta^2 D^3)
\]
satisfies $\wt(\hat{v})=\wt(\hat{v}-v_1)=4=d(\hat{v},\C)$ for the codeword
\[
  v_1=(1+\beta D+\beta^5 D^4+\beta^2D^5,\,1+\beta^3D^2+\beta^2D^3+\beta^2D^5)\in\C.
\]
Hence $\hat{v}$ is equally close to~$v_1$ and the zero codeword, but only
$\hat{v}-v_1$ satisfies the error condition\eqr{e-slide.error}.
Therefore, the decoding algorithm will decode~$\hat{v}$ into the codeword~$v_1$.
\end{rem}

The main step of the algorithm in Theorem~\ref{T-dec.alg} is, of course, the
determination of $\eta_0$ from the syndrome vector in Step~3.
For codes with small parameters this can easily be achieved by simply
checking (in a smart way) all linear combinations of at most $\frac{M+1}{2}$
columns of $H^c_M$. But for larger codes this is unsatisfactory and one would
like to know an algebraic computation of $\eta_0$.
Unfortunately, thus far we cannot offer such an algebraic decoding.
It will certainly depend on an algebraic construction of superregular
matrices along with their algebraic properties.

We close this paper with the following criterion which, in the affirmative
case, speeds up Step~3.
It makes use of the systematic sliding parity check matrix of
the code, see\eqr{e-hatH1}, which can be used just as well in the decoding
algorithm.
Notice that there are~$n$ different
systematic $M$-th sliding parity check matrices for an $(n,n-1,\delta)$
code. Therefore, the following criterion can be tested~$n$ times and
hopefully leads to an immediate decision on $\eta_0$ at least ones.

\begin{prop}\label{P-S.eta0}
Let $\C\subseteq\FDlaurent^n$ be a strongly MDS $(n,n-1,\delta)$-code with
systematic $M$-th sliding parity check matrix~$\hat{H}$ as in\eqr{e-hatH1}.
Let $\hat{S}=(\hat{S}_0,\ldots,\hat{S}_M)\in\F^{M+1}$ be such that
$\hat{S}=\hat{\eta}\hat{H}\T$ for some
\[
  \hat{\eta}=(e_0,\ldots,e_M,f_0,\ldots,f_M)\in\F^{(M+1)+(M+1)(n-1)}\text{ and }
  \wt(\hat{\eta})\leq\frac{M+1}{2}.
\]
If $\wt(\hat{S})\leq\lceil\frac{M+1}{2}\rceil$, then $f_0=0$ and $e_0=S_0$.
\end{prop}

\begin{proof}
The assumptions $\wt(\hat{S})\leq\lceil\frac{M+1}{2}\rceil$ and
$\wt(\hat{\eta})\leq\frac{M+1}{2}$ imply that there exists a linear
combination of at most $M+1$ columns of~$\hat{H}$ giving the zero vector.
But then Theorem~\ref{T-strongMDS} yields $f_0=0$ and $e_0=S_0$.
\end{proof}

\section{Conclusion}

In this paper we introduced two new classes of convolutional
codes called strongly MDS convolutional codes and codes having
maximum distance profile. Strongly MDS convolutional codes have
the property that the generalized Singleton bound is attained at
the earliest possible column distance. Codes with maximum
distance profile have a maximal possible increase of the column
distances.

It is shown that strongly MDS convolutional codes perform
excellent under feedback decoding. The number of errors which can
be iteratively decoded per time interval lets these codes compare
with MDS linear block codes having a considerable block
length. At this point the feedback decoding algorithm we
presented is not powerful enough to practically decode strongly
MDS convolutional codes when the degree $\delta$ and the block
length $n$ are too large for the syndrome decoding step, see
Theorem~\ref{T-dec.alg} and Proposition~\ref{P-S.eta0} for
details. It will be a matter of future research to construct
strongly MDS convolutional codes which come equipped with an
algebraic structure and an efficient decoding algorithm
comparable to the situation of Reed-Solomon block codes. The
class of cyclic convolutional codes~\cite{gl02p,gl02u,pi76,ro79}
might hold some promise here.

{}From an applications point of view strongly MDS convolutional
codes are particularly suited in situations where codes over
large alphabets are required and in situations where algebraic
decoding is desirable. Hadjicostis~\cite{ha02,ha03a} has recently
demonstrated that convolutional codes over large alphabets are
very desirable in areas of process control via linear finite
state machines where large numbers of non-concurrent errors
should be detected and corrected. It seems that strongly MDS
convolutional codes have potential for such applications.

\section*{Appendix}
\renewcommand{\theequation}{A.\arabic{equation}}
\setcounter{equation}{0}

We will prove that the proper minors of the matrix $T_n$ given in Example~\ref{E-superreg2}(2)
are all positive.
In order to do so consider the matrix
\[
X=\begin{bmatrix}
  1&&&&&\\
  1&1&&&&\\

  &1&1&&&\\
  &&\ddots&\ddots &&\\
  &&&1&1&\\
  &&&&1&1\\
  \end{bmatrix}\in\Z^{n\times n}
\]
and notice that for all $k\in\{1,\ldots,n-1\}$ we have
\begin{equation}\label{e-Xk}
X^k=\begin{bmatrix}
   1&&&&&&&&\\
   {\binom{k}{1}}&1&&&&&&&\\
   {\binom{k}{2}}&{k \choose 1}&1&&&&&&\\
   \vdots&\ddots &\ddots&\ddots &&&&&\\
   \vdots& &\ddots&\ddots &\ddots&&&&\\
   1&\ldots  &\ldots&{k \choose 2}&{k \choose 1}&1&&&\\
    &1&\ldots&\ldots& {k \choose 2}&{k \choose 1}&1&&\\
   &&\ddots&&&\ddots&\ddots&\ddots&\\
   &&&1&\ldots&\ldots&{k \choose 2}&{k \choose 1}&1\\
  \end{bmatrix}.
\end{equation}
In particular, $X^{n-1}=T_n$.
Therefore, the positivity of the proper minors is a consequence of the following theorem.

\noindent
{\bf Theorem A}
{\sl Let $k\in\{1,\ldots,n-1\}$ and $1\leq i_1<\ldots<i_r\leq n,\,1\leq j_1<\ldots<j_r\leq n$ and
define $\hat{X}:=(X^k)^{i_1,\ldots,i_r}_{j_1,\ldots,j_r}$. Then
$\det\hat{X}\geq0$ and
\[
   \det\hat{X}>0\Longleftrightarrow j_l\in\{i_l,i_l-1,\ldots,i_l-k\}
   \text{ for all }l=1,\ldots,r.
\]
}

\begin{proof}
1) We first show that
\begin{equation}\label{e-zerominors}
  j_l\not\in\{i_l,i_l-1,\ldots,i_l-k\}\text{ for some }l
  \Longrightarrow\det\hat{X}=0.
\end{equation}
To this end notice that
\[
    X_{ij}=0\text{ for } j>i \text{ or }j<i-k
\]
and thus
\[
   \hat{X}_{ef}=X_{i_e j_f}=0\text{ for } j_f>i_e \text{ or }
         j_f<i_e-k.
\]
Assume now $j_l>i_l$ for some~$l$. Then for all $e\leq l$ and $f\geq l$ we have
$j_f\geq j_l>i_l\geq i_e$ and thus $\hat{X}_{ef}=0$.
Hence the first~$l$ rows of~$\hat{X}$ have at most rank~$l-1$ and thus
$\det\hat{X}=0$.
Similarly, if $j_l<i_l-k$ for some~$l$, then we have $\hat{X}_{ef}=0$ for all
$e\geq l$ and $f\leq l$ and the first~$l$ columns of~$\hat{X}$ have at most
rank~$l-1$.

2) It remains to prove the implication ``$\Longleftarrow$'' of the equivalence given
in the theorem.
\\
We begin with proving the statement for $k=1$, i.~e. for the matrix $X$.
In order to do so, we proceed by induction on~$r$.
For $r=1$ we have to consider the submatrices $X^{i_1}_{i_1}$ and $X^{i_1}_{i_1-1}$.
They all trivially have determinant~$1$.
Now let $r>1$. We suppose the statement is true for all $(r-1)\times (r-1)$ proper
submatrices with the according restriction on the indices and we have to show that the assertion
is also true for $\hat{X}=X^{i_1,\ldots,i_r}_{j_1,\ldots,j_r}$ where $j_l\in\{i_l,i_l-1\}$ for
all~$l$.
Notice that the first column of~$\hat{X}$ has either
one or two nonzero entries and they are equal to~$1$.
If the first column of $\hat{X}$ has one~$1$ only, then it is on the
first row. Applying cofactor expansion along that column we
obtain
\begin{equation}\label{e-expansion}
  \det\hat{X} =1\cdot \det~X^{i_2,\ldots,i_r}_{j_2,\ldots,j_r}.
\end{equation}
The $(r-1)\times (r-1)$-submatrix satisfies
$j_{l}\in\{i_{l},i_{l}-1\} $ for all $l=2,\ldots,r$ and hence by induction has
positive determinant.
This proves $\det\hat{X}>0$ in this case.
If the first column of $X^{i_1,\ldots,i_r}_{j_1,\ldots,j_r}$
has two entries equal to~$1$, then
they are necessarily on the first two rows, thus $i_2=i_1+1$ and $j_1=i_1$.
Since $j_{2}\in \{ i_{2},i_{2}-1\}=\{i_{1}+1,i_{1}\}$ and $j_2> j_1$,
we can only have $j_2=i_1+1.$
Then the first row will have only one nonzero entry equal to~$1$
on the first position, and applying cofactor expansion
along that row, we obtain again\eqr{e-expansion} and thus $\det\hat{X}>0$.
\\
We now proceed by induction on~$k$ in order to prove the desired result for $X^k$
where $k>1$. Assume  $X^{k-1}$ has the stated property.
Using $X^k=X\cdot X^{k-1}$ and the Cauchy-Binet formula for minors
we obtain
\[
  \det\hat{X}=
  \sum_{1\leq s_1<\ldots< s_r\leq n,\atop s_l\in \{i_l,i_{l-1}\}\cap \{j_l,j_{l}+1,\ldots,
    j_l+{k-1}\}}
  \det X^{i_1,\ldots,i_r}_{s_1,\ldots,s_r}\cdot
 \det(X^{k-1})^{s_1,\ldots,s_r}_{j_1,\ldots,j_r}.
\]
Due to part~1) of the proof the sum indeed expands only over the given indices.
By induction all nonsingular submatrices of both matrices~$X$ and~$X^{k-1}$ have
positive determinant, hence if there are any nonzero terms in the sum, it is
necessarily positive.
Therefore, the only thing left to be proven is that there
is a nonzero term in the above sum. But all products of the form
$\det X^{i_1,\ldots,i_r}_{i_1,\ldots,i_r}\cdot
 \det(X^{k-1})^{i_1,\ldots,i_r}_{j_1,\ldots,j_r}$ with
$j_{l}\in\{ i_{l},i_{l}-1,i_{l}-2,\ldots,  i_{l}-(k-1)\}$ for all~$l$
are nonzero. Thus $\det\hat{X}>0$ and the proof is complete.
\end{proof}\bigskip

{\sc Proof of Theorem~\ref{T-str.MDSn}:}
\underline{Step 1:}
We will show the existence of a systematic sliding parity check matrix~$\hat{H}$
as in\eqr{e-hatH1} satisfying part~(b) of Theorem~\ref{T-strongMDS}.
This can be accomplished as follows.
Let $\tau:=(M+1)(n-1)$ and pick a
$\tau\times\tau$-superregular matrix in Toeplitz form, say
\[
    T:=\begin{bmatrix}t_1    & 0      & \cdots & 0 \\
                    t_2    & t_1    & \ddots & \vdots  \\
                    \vdots & \ddots &\ddots & 0\\
                    t_{\tau}    & \cdots &   t_2  & t_1\end{bmatrix}
     =\begin{bmatrix} T^1\\ T^2\!\!\phantom{\vdots}\\ \vdots\\ T^{\tau}\end{bmatrix}
     =\left[T_1,\ldots, T_{\tau}\right].
\]
Theorem~\ref{T-exist} guarantees the existence of such a matrix over a suitably large
field of characteristic~$p$.
Now define
\[
  \hat{H}=\left[\!\begin{array}{cccc|c}
      1&      &      & &  T^{n-1}  \\
       &1     &      & &  T^{2(n-1)} \\
       &      &\ddots& & \vdots \\
       &      &      &1&  T^{(M+1)(n-1)}
  \end{array}\!\right]\in\F^{(M+1)\times(M+1)n}.
\]
Notice that by construction~$\hat{H}$ has the form as in\eqr{e-hatH1}.
We will prove by contradiction that this matrix satisfies part~(b) of
Theorem~\ref{T-strongMDS}.
In order to do so, write
$\hat{H}=[e_1,\ldots,e_{M+1},\hat{T}_1,\ldots,\hat{T}_{\tau}]$ and
assume that $i\leq n-1$ is the smallest index such that~$\hat{T}_i$ is in the
span of~$M$ other columns of~$\hat{H}$. Hence these other columns do not involve
$\hat{T}_1,\ldots,\hat{T}_{i-1}$.
This implies that there is a linear combination of~$M+1$ columns of the
matrix $[I_{\tau},T]$ with a nonzero coefficient for the column~$T_i$ and having a zero entry
at the positions $1,2,\ldots,i-1,n-1,2(n-1),\ldots,(M+1)(n-1)$.
Since $i\leq n-1$, these positions are indeed different and thus the
weight of this linear combination is at most~$\tau-i+1-(M+1)$.
Consider now the matrix
\[
   Y:=[I_{\tau-i+1}\,|\, \tilde{T}_i,\tilde{T}_{i+1},\ldots,\tilde{T}_{\tau}]
     :=\left[\!\begin{array}{ccc|ccc}
      1&      & &t_1&  & \\
       &\ddots& &\vdots&\ddots& \\
       &      &1&t_{\tau-i+1}&\cdots&t_1
  \end{array}\!\right]
   \in\F^{(\tau-i+1)\times2(\tau-i+1)},
\]
where $\tilde{T}_j$ denotes vector of the last $\tau-i+1$ entries
of~$T_j$.
Notice that superregularity of~$T$ implies superregularity of the matrix
$[\tilde{T}_i,\ldots,\tilde{T}_{\tau}]$.
The linear combination of $M+1$ columns of~$[I_{\tau},T]$ above now reads
as a linear combination of~$M+1$ columns of~$Y$ with a nonzero coefficient for
the column $\tilde{T}_i$ and having weight at most $\tau-i+1-(M+1)$.
Hence picking a suitable set of (at most) $\tau-i+1-(M+1)$ standard basis vectors,
we obtain that the column~$\tilde{T}_i$ is in the span of $\tau-1$ other columns of~$Y$.
But this is a contradiction to Theorem~\ref{T-superreg}(d).

\underline{Step 2:}
Having constructed a matrix~$\hat{H}$ as in\eqr{e-hatH1} with the corresponding
column condition, we now establish the existence of an $(n,n-1)$-code
having~$\hat{H}$ as its $M$-th systematic sliding parity check
matrix.
In order to simplify notation write
\begin{eqnarray}   \label{e-hatH3}
  \hat{H}&=&\left[\!\begin{array}{cccc|cccc}
                   1&  &      &  &h_0    & 0      & \cdots & 0 \\
                    & 1&      &  &h_1    & h_0   & \ddots & \vdots \\
                    &  &\ddots&  & \vdots & \ddots &\ddots & 0\\
                    &  &      &1 &h_M& \cdots &
                    h_1&h_0\end{array}\!\right]  \nonumber \\
        &=&[e_1,\ldots,e_{M+1},H_{12},\ldots,H_{1n},\ldots,H_{M+1,2},\ldots,H_{M+1,n}]
\end{eqnarray}
where $h_i=(h_{i2},\ldots,h_{in})\in\F^{n-1}$.
We have to find polynomials
\begin{equation}\label{e-ab}
   a=1+\sum_{i=1}^{\delta}a_iD^i\in\F[D],\quad
   b=\sum_{i=0}^{\delta}b_iD^i\in\F[D]^{n-1}
\end{equation}
such that
\begin{equation}\label{e-Laurent3}
    \frac{b}{a}=\sum_{j=0}^Mh_{j}D^j+\text{ higher powers}
\end{equation}
(see Theorem~\ref{T-strongMDS}).
Recall that $M=\lfloor\frac{\delta}{n-1}\rfloor+\delta$.
If $M=\delta$ (i.~e. $\delta<n-1$), we may simply take $a=1$ and
$b=\sum_{i=0}^{\delta}h_iD^i$.
Now let us assume $M>\delta$.
Comparing like powers of~$D$ in\eqr{e-Laurent3} shows that the above requires in
particular
\begin{equation}\label{e-coeff}
    0=h_l+a_1h_{l-1}+a_2h_{l-2}+\ldots+a_{\delta}h_{l-\delta}
    \text{ for any }l>\delta
\end{equation}
for suitable $h_{M+1}, h_{M+2},\ldots\in\F^{n-1}$.
For $l=\delta+1,\ldots,M$ these equations read as
\begin{equation}\label{e-rowsp}
    (a_{\delta},\ldots,a_1)
    \begin{bmatrix} h_{M-\delta}&h_{M-\delta-1}&\cdots& h_1\\
                   h_{M-\delta+1}&h_{M-\delta}&       &h_2\\
                    \vdots &\vdots&   &\vdots \\
                    h_{M-1}&h_{M-2}&\cdots&h_{\delta}\end{bmatrix}
    =-(h_M,\ldots,h_{\delta+1}).
\end{equation}
Notice that $h_1,\ldots,h_M$ are given data.
If we can find a solution $(a_{\delta},\ldots,a_1)$ of\eqr{e-rowsp},
then\eqr{e-coeff} can be established for all $l\geq M+1$ by choosing $h_l$ suitably.
Thereafter, the vector polynomial $b\in\F[D]^{n-1}$ can be computed by
equating the coefficients of $D^0,\ldots,D^{\delta}$ in the equation
$b=(\sum_{i=0}^{\infty}h_iD^i)(1+\sum_{i=1}^{\delta}a_i D^i)$.
Hence it remains to consider\eqr{e-rowsp}.
This equation is solvable if
\[
   \rank\H=(n-1)(M-\delta),\ \text{ where }
   \H:=\begin{bmatrix} h_{M-\delta}&\cdots& h_1\\ \vdots &   &\vdots \\
                    h_{M-1}&\cdots&h_{\delta}\end{bmatrix}
   \in\F^{\delta\times(n-1)(M-\delta)}.
\]
Notice that $\rho:=(n-1)(M-\delta)=(n-1)\lfloor\frac{\delta}{n-1}\rfloor\leq\delta$.
We proceed by contradiction and assume $\rank\H<\rho$.
Then there is a column of~$\H$ that is a linear combination of the other $\rho-1$
columns.
Since~$\H$ is a submatrix of~$\hat{H}$ (see\eqr{e-hatH3}) and because of the specific
structure of $\H$, this yields that a
column~$H_{1j},\;j=2,\ldots,n$, is a linear combination of $\rho-1+M+1-\delta$
other columns of~$\hat{H}$.
But
\[
  M+\rho-\delta=\Big\lfloor\frac{\delta}{n-1}\Big\rfloor+
  \delta+(n-1)\Big\lfloor\frac{\delta}{n-1}\Big\rfloor-\delta
  \leq\Big\lfloor\frac{\delta}{n-1}\Big\rfloor+\delta=M,
\]
and thus we arrive at a contradiction to the column property of~$\hat{H}$ (see~(b) of
Theorem~\ref{T-strongMDS}).
Hence\eqr{e-rowsp} is solvable and the existence of~$a$ and~$b$ as in\eqr{e-ab}
and\eqr{e-Laurent3} is established.

\underline{Step 3:}
Put $H=[a,b^{(1)},\ldots,b^{(n-1)}]$, where~$a\in\F[D]$ and
$b=:(b^{(1)},\ldots,b^{(n-1)})\in\F[D]^{n-1}$ are constructed as in Step~2).
Moreover, let $\C=\{v\in\FDlaurent^n\mid v H\T=0\}$.
It remains to show that~$\C$ has degree~$\delta$,
which amounts to showing that $a,\,b^{(1)},\ldots,b^{(n-1)}$ are coprime and
\begin{equation}\label{e-degree}
  \max\{\deg a ,\,\deg b^{(1)},\ldots,\,\deg b^{(n-1)}\}=\delta.
\end{equation}
Coprimeness can be assumed without loss of generality since division by a
common factor would lead to another solution of\eqr{e-ab} and\eqr{e-Laurent3}.
Hence~$H$ is basic.
By construction and Theorem~\ref{T-strongMDS} the $M$-th column distance of $\C$ is
given by
$d^c_M=\lfloor\frac{\delta}{n-1}\rfloor+\delta+2$.
Since this number is strictly bigger than the generalized Singleton bound of any
$(n,n-1,\hat{\delta})$-code, where $\hat{\delta}<\delta$, Equation\eqr{e-degree}
follows immediately.

Thus $\C$ is a strongly MDS $(n,n-1,\delta)$-code and the proof is
complete.
\hfill$\Box$\par


\end{document}